\documentclass[journal,twoside]{IEEEtran}

\usepackage{amsfonts,amsmath,amssymb,mathrsfs,float,latexsym,amstext}
\usepackage{amsthm}
\usepackage{graphics,graphicx}
\RequirePackage[OT1]{fontenc}
\usepackage{cite}

\usepackage[colorlinks=true,linkcolor=black,anchorcolor=black,citecolor=black,filecolor=black,menucolor=black,runcolor=black,urlcolor=black]{hyperref}

\newcommand{\ignore}[1]{} %%% {} empty inside

\allowdisplaybreaks %% allows to break long multiline formulas (use with align, not with aligned!)

\theoremstyle{plain}
\newtheorem{theorem}{Theorem}%[section]
\newtheorem{proposition}{Proposition}
\newtheorem{lemma}{Lemma}%[section]

\theoremstyle{definition}

\newtheorem{remark}{Remark}%[section] 
\font\msbmx=msbm10                   % \emptyset should be changed to \varnothing
\font\msbmvii=msbm7                  % in the paper and this will give the proper symbol
\font\msbmv=msbm5
\def\varnothing{\mathchoice{\mbox{\msbmx\char'077}}%
{\mbox{\msbmx\char'077}}{\mbox{\msbmvii\char'077}}{\mbox{\msbmv\char'077}}}%
\def\Var{\mathop{\rm Var}\nolimits}%   a math operator.

\newcommand{\mb}[1]{\mathbf{#1}} %% use {\mb{symb}} provides mathbold roman
\newcommand{\Cb}{{\mb{C}}}

\newcommand{\Xb}{{\mb{X}}}

\def\C{{\mb{C}}}
\def\CP{{\mb{CP}}}

\newcommand{\mbs}[1]{\boldsymbol{#1}} %% use {\mbs{symb}} provides mathbold italic, better for greek symbols

\newcommand{\alphab}{{\mbs{\alpha}}}
\newcommand{\betab}{{\mbs{\beta}}}

\newcommand{\mc}[1]{\mathcal{#1}} %% use \mc{symb}
\newcommand{\Rc}{{\mc{R}}}
\newcommand{\Rca}{{\bar{\mc{R}}}}

\newcommand{\Nc}{{\mc{N}}}
\newcommand{\Fc}{{\mc{F}}}

\newcommand{\Ac}{{\mc{A}}}

\newcommand{\Xc}{{\mc{X}}}
\newcommand{\Bc}{{\mc{B}}}

\newcommand{\Lc}{{\mc{L}}}
\newcommand{\Kc}{{\mc{K}}}

\newcommand{\wtX}{\widetilde{X}}
\newcommand{\wtS}{\widetilde{S}}

\usepackage{mathrsfs}
 \newcommand{\PFA}{\mathsf{PFA}}% prob of false alarm
 \newcommand{\PMI}{\mathsf{PMI}}% prob of false alarm
\newcommand{\Pb}{{\mathsf{P}}} %probability
\newcommand{\wtPb}{\widetilde{\Pb}}

\newcommand{\Eb}{{\mathsf{E}}}
\newcommand{\wtE}{\widetilde{\Eb}}
\newcommand{\Hyp}{{\mathsf{H}}} %  hypothesis
\newcommand{\mrm}[1]{\mathrm{#1}}

\newcommand{\D}{{\mrm{d}}}

\newcommand{\mbb}[1]{\mathbb{#1}} %use \mbb{symb}
\def\One{\mathchoice{\rm 1\mskip-4.2mu l}{\rm 1\mskip-4.2mu l}
{\rm 1\mskip-4.6mu l}{\rm 1\mskip-5.2mu l}}
\newcommand\Ind[1]{{\One_{\{#1\}}}} % indicator -- use \Ind{A}
\newcommand{\Zbb}{\mbb{Z}} %% discrete real line

\newcommand{\class}{{\mbb{C}}}

\newcommand{\xra}{\xrightarrow} %% use \xla{\text{prob}} matches the length of arrow

\newcommand{\abs}[1]{\left\vert#1\right\vert}
\newcommand{\set}[1]{\left\{#1\right\}}

\newcommand{\brc}[1]{\left(#1\right)}
\newcommand{\brcs}[1]{\left[#1\right]}

\renewcommand{\le}{\leqslant} % AMS le ge
\renewcommand{\ge}{\geqslant}

\usepackage{color}
\begin{document}
%%%%%%%%%%%%%%%%

%%%%
\title{An Asymptotic Theory of Joint Sequential Changepoint Detection and Identification for General Stochastic Models
\thanks{The work was supported in part by the Russian Science Foundation under the grant 18-19-00452
at the Moscow Institute of Physics and Technology.}}

%%%%
\author{Alexander~G.~ Tartakovsky, \IEEEmembership{Senior~Member,~IEEE}  
 \thanks{ A. G.  Tartakovsky is a Deputy Head of the Space informatics Laboratory at the Moscow Institute of Physics and Technology, Russia
 and President of AGT StatConsult, Los Angeles, California, USA; e-mail: agt@phystech.edu}
\thanks{Manuscript received  March 23, 2020; revised October 26, 2020; accepted February 22, 2021.}
\thanks{Copyright (c) 2020 IEEE. Personal use of this material is permitted.  However, permission to use this material for any other purposes must be obtained from the 
IEEE by sending a request to pubs-permissions@ieee.org.}
}

%%%
\markboth{IEEE Transactions on Information Theory,~Vol.~~, No.~~, ~~~2021}%
{Tartakovsky: Joint Sequential Changepoint Detection and Identification }

\maketitle

%%%
\begin{abstract}
The paper addresses a joint sequential changepoint detection and identification/isolation problem for a general stochastic model, assuming that the observed data may be  dependent and 
non-identically distributed, the prior distribution of the change point is arbitrary, and the post-change hypotheses are composite. The developed detection--identification theory generalizes  
the changepoint detection theory developed by Tartakovsky (2019) to the case of multiple composite post-change hypotheses when one has not only to detect a change as quickly as possible but also to identify 
(or isolate) the true post-change distribution.  We propose a multi-hypothesis change detection--identification rule and show that it is nearly optimal, minimizing moments of the delay to detection as 
the probability of a false alarm and the probabilities of misidentification go to zero.   
\end{abstract}

%%%

%%%%%%%%%%%%%%%%%%%%%%%%%%%%%%%%%%%%%%%%%%%%%%
%_____________________________________
\section{Introduction} \label{sec:intro}

\IEEEPARstart{I}{n} many applications, one needs not only to detect an abrupt change as quickly as possible but also to provide a detailed diagnosis of the occurred change -- to determine which type of change is in effect. 
For example, the problem of detection and diagnosis is important for rapid detection and isolation of intrusions in large-scale distributed computer networks, target detection with radar, sonar and optical sensors
in a cluttered environment, detecting terrorists' malicious activity, fault detection and isolation in dynamic systems and networks, and integrity monitoring of navigation systems, to name a few (see \cite[Ch 10]{TNB_book2014}
for an overview and references).
In other words, there are several kinds of changes that can be associated with several different post-change distributions and the goal is to detect the change and to identify which distribution corresponds to the change. 
As a result, the problem of changepoint detection and diagnosis is a generalization of the quickest change detection problem \cite{ShiryaevTPA63, ShiryaevBook78, lorden-ams71,PollakAS85,TNB_book2014} 
to the case of $N \ge 2$ post-change hypotheses, 
and it can be formulated as a joint change detection and identification problem. In the literature, this problem is usually called {\em change detection and isolation}. The detection--isolation problem has been considered in
both Bayesian and minimax settings. In 1995, Nikiforov~\cite{NikiforovIEEEIT95} suggested a minimax approach to the change detection--isolation problem and showed that the multihypothesis version of the CUSUM rule is
asymptotically optimal when the average run length (ARL) to a false alarm and the mean time to false isolation become large. Several versions of the multihypothesis CUSUM-type and SR-type procedures, which have
minimax optimality properties in the classes of rules with constraints imposed on the ARL to a false alarm and conditional probabilities of false isolation, are proposed by Nikiforov~\cite{NikiforovIEEEIT00, NikiforovIEEEIT03}
and Tartakovsky~\cite{Tartakovsky-SQA08b}. These rules asymptotically minimize maximal expected delays to detection and isolation as the ARL to a false alarm is large and the probabilities of wrong
isolations are small. Dayanik~et al.~\cite{dayaniketal-AOR13} proposed an asymptotically optimal Bayesian detection--isolation rule assuming that the prior distribution of the change point is geometric. 
In all these papers, the optimality results were restricted to the case of independent and identically distributed (i.i.d.) observations (in pre- and post-change modes with different distributions) and simple 
post-change hypotheses. In many practical applications, the i.i.d.\ assumption is too restrictive. The observations may be
either non-identically distributed or dependent or both, i.e., non-i.i.d.  Also, in a variety of applications, a pre-change distribution is known but the post-change distribution is rarely known completely. 
A more realistic situation is parametric uncertainty when the parameter  of the post-change distribution is unknown since a putative parameter value is rarely representative. 
Lai~\cite{LaiIEEE00} provided a certain generalization for the non-i.i.d.\ case and composite hypotheses for a specific loss function. 
See Chapter~10 in Tartakovsky et al.~\cite{TNB_book2014} for a detailed overview.

One of the most challenging and important versions of the change detection--isolation problem is the multidecision and multistream detection
problem when it is necessary not only to detect a change as soon as possible but also  to identify the streams where the change happens with given probabilities of misidentification.
Specifically, there are $N$ data streams and the change occurs in some of them at an unknown point in time.  It is necessary to detect
the change in distribution as soon as possible and indicate which streams are ``corrupted.'' Both the rates of false alarms and misidentification should be controlled by given (usually low) levels. In the following, we will refer 
to this problem as the  {\em Multistream Sequential Change Detection--Identification} problem. 

In this paper, we address a simplified multistream detection--identification scenario where change can occur only in a single stream and we need to determine in which stream. We assume that the
observations in streams can have a very general structure, i.e., can be dependent and non-identically distributed. We focus on a semi-Bayesian setting assuming that the change point is random possessing the 
(prior) distribution.
However, we do not suppose that there is a prior distribution on post-change hypotheses. We generalize the asymptotic Bayesian theory developed by Tartakovsky~\cite{TartakovskyIEEEIT2019} for a single 
post-change hypothesis (for a single stream). Specifically, we show that under certain 
conditions (related to the law of large numbers for the log-likelihood processes) the proposed multihypothesis detection--identification rule asymptotically minimizes the trade-off between positive moments of the detection
delay and the false alarm/misclassification rates expressed via the weighted probabilities of false alarm and false identification. The key assumption in the general asymptotic theory is the 
stability property of the log-likelihood ratio processes in streams between the  ``change'' and ``no-change'' hypotheses, which can be formulated in terms of the law of large numbers and rates of convergence
of the properly normalized log-likelihood ratios and their adaptive versions in the vicinity of the true parameter values. 

The rest of the paper is organized as follows.  In Section~\ref{sec:Model}, we describe the general stochastic model, which is treated in the paper. In Section~\ref{sec:Rule}, we introduce the mixture-based 
change detection--identification rule.  In Section~\ref{sec:Ch8Problems}, we formulate the asymptotic optimization problems in the class of changepoint detection--identification rules with the constraint imposed 
on  the probabilities of false alarm and wrong identification. In Section~\ref{sec:PFAandPMI}, we obtain upper bounds on the probabilities of false alarms and misidentification as functions of thresholds. 
In Section~\ref{sec:LBdetisol}, we  derive asymptotic lower bounds for moments of the detection delay in the class of rules with given probabilities of false alarms and misidentification, 
and  in Section~\ref{sec:Asoptdetisol}, we prove asymptotic optimality of the proposed mixture detection--identification rule as the probabilities of false alarm and misidentification go to zero.
In Section~\ref{sec:Ex}, we consider an example that illustrates general results. Section~\ref{sec:Remarks} concludes.

%%____________________________________________________________________
\section{The General Stochastic Model}\label{sec:Model}

Suppose there are $N$ independent data streams $\{X_n(i)\}_{n \ge 1}$, $i=1,\dots,N$, observed sequentially in time subject to a change at an unknown time $\nu\in\{0,1, 2, \dots\}$, 
so that $X_1(i),\dots,X_\nu(i)$ are generated by one stochastic model and $X_{\nu+1}(i),  X_{\nu+2}(i), \dots$ by another model when the change occurs in the $i$th stream. 
We will assume that the change in distributions may happen only in one stream and it is not known which stream is affected, i.e., we are interested in
a ``multisample slippage'' changepoint model (given $\nu$ and that the $i$th stream is affected with the parameter $\theta_i$) for which joint density $p(\Xb^n | \Hyp_{\nu,i}, \theta_i)$ of the data 
$\Xb^n= (\Xb^n(1),\dots,\Xb^n(N))$, $\Xb^n(i)= (X_1(i),\dots,X_n(i))$ observed up to time $n$ is of the form 
\begin{equation}\label{modeldetisol}
\begin{split}
& p(\Xb^n | \Hyp_{\nu,i}, \theta_i)  =  p(\Xb^n| \Hyp_\infty)
\\
& = \prod_{t=1}^n \prod_{\ell=1}^N g_\ell(X_t(\ell)|\Xb^{t-1}(\ell)) ~~\text{for}~~ \nu \ge n ,
\\
& p(\Xb^n | \Hyp_{\nu,i},\theta_i)  =  \prod_{t=1}^{\nu} g_i(X_t(i)|\Xb^{t-1}(i))  \times
\\
& \prod_{t=\nu+1}^{n}  f_{i,\theta_i}(X_t(i)|\Xb^{t-1}(i)) \times 
\\
& \prod_{j\in \Nc\setminus \{i\}}   \prod_{t=1}^n g_j(X_t(j)|\Xb^{t-1}(j)) ~~ \text{for}~~ \nu < n,
\end{split}
\end{equation}
where $\Hyp_{\nu,i}$ denotes the hypothesis that the change occurs at time $\nu$ in the stream $i$, $g_i(X_t(i)|\Xb^{t-1}(i))$ and $f_{i,\theta_i}(X_t(i)|\Xb^{t-1}(i))$ are conditional pre- and post-change densities in the $i$th data stream, respectively (with respect to some sigma-finite measure), and $\Nc=\{1,2,\dots,N\}$. In other words,  all components $X_t(\ell)$, $\ell\in \Nc$,
have conditional densities $g_\ell(X_t(\ell)| \Xb^{t-1}(\ell))$ before the change occurs and $X_t(i)$  has conditional density $f_{i,\theta_i}(X_t(i)| \Xb^{t-1}(i))$ after the change occurs in the $i$th stream and the rest of the components 
$X_t(j)$, $j\in \Nc\setminus\{i\}$  have conditional densities  $g_j(X_t(j)| \Xb^{t-1}(j))$. The parameters $\theta_i\in\Theta_i$, $i=1,\dots,N$ of the post-change distributions are unknown.  The event $\nu=\infty$ and the
corresponding hypothesis $\Hyp_\infty: \nu=\infty$ mean that there never is a change. Notice that the model \eqref{modeldetisol} implies that $X_{\nu+1}(i)$ is the first post-change observation under hypothesis $\Hyp_{\nu,i}$.

Regarding the change point $\nu$ we assume that it is a random variable independent of the observations with prior distribution 
$\pi_k=\Pb(\nu=k)$, $k=0,1,2,\dots$ with $\pi_k >0$ for $k\in\{0,1,2, \dots\}=\Zbb_+$. We will also assume that a change point may take negative values, which means that the change has 
occurred by the time the observations became available. However, 
the detailed structure of the distribution $\Pb(\nu=k)$ for $k=-1,-2,\dots$ is not important. The only value which matters is the total probability $q=\Pb(\nu \le -1)$ of the change being in 
effect before the observations become available,  so we set $\Pb(\nu \le -1)= \Pb(\nu=-1)=\pi_{-1}$, $ \pi_{-1} \in [0,1)$.  Therefore, in what follows we assume that $\nu\in \{-1,0,1,\dots\}=\Zbb$ and the prior distribution of
the change point is defined on $\Zbb$.

%%____________________________________________________________________
\section{The Detection--Identification Rule}\label{sec:Rule}

A changepoint detection--identification rule is a pair $\delta=(d, T)$, where $T$ is a stopping time (with respect to the filtration $\{\Fc_n=\sigma(\Xb^n)\}_{n \in \Zbb_+}$) associated with the time of alarm on change and 
$d=d_T \in\Nc$ is a decision on which stream is affected (or which post-change distribution is true) which is made at time $T$.

It follows from \eqref{modeldetisol} that for an assumed value of the change point $\nu = k$,  stream $i \in \Nc$, and the post-change parameter value in the $i$th stream $\theta_i\in\Theta_i$, the 
likelihood ratio (LR)
$LR_{i,\theta_i}(k,n)= p(\Xb^n | \Hyp_{k,i},\theta_i)/p(\Xb^n| \Hyp_\infty)$ between the hypotheses $\Hyp_{k,i}$ and $\Hyp_\infty$ for observations accumulated by the time $n$ has the form
\begin{equation}\label{LRdeteisol}
 LR_{i,\theta_i}(k,n)= \prod_{t=k+1}^{n}  \Lc_{i,\theta_i}(t), \quad i\in \Nc, ~~ n > k
\end{equation}
($k =-1,0,1, \dots$), where 
\[
\Lc_{i,\theta_i}(t)=f_{i,\theta_i}(X_t(i)|\Xb^{t-1}(i))/g_i(X_t(i)|\Xb^{t-1}(i)).
\]
We suppose that $\Lc_{i,\theta_i}(0)=1$, so that $LR_{i,\theta_i}(-1,n)=LR_{i,\theta_i}(0,n)$. 
Define the average (over the prior $\pi_k$) LR statistics
\begin{equation}\label{LRpiskn} 
\Lambda_{i,\theta_i}^\pi(n)  =   \sum_{k=-1}^{n-1} \pi_k LR_{i,\theta_i}(k, n) , \quad  i \in \Nc .
\end{equation}

Let $W_i(\theta_i)$, $\int_{\Theta_i} {\rm d}  W_i(\theta_i) =1$, $i\in \Nc$ be mixing measures  and, for $k<n$ and $i\in\Nc$,
define the LR-mixtures  
\begin{equation}\label{LRmixtureskn} 
LR_{i, W}(k,n)  =  \int_{\Theta_i} LR_{i,\theta_i}(k, n)  \, \mrm{d} W_i(\theta_i),  
\end{equation}
and the statistics 
\begin{equation}\label{LRmixturespi}
\Lambda_{i, W}^\pi(n) =
\begin{cases}
\sum_{k=-1}^{n-1} \pi_k LR_{i,W}(k,n), ~~ & i\in \Nc
\\
\Pb(\nu \ge n) ~~ & i=0
\end{cases} ;
\end{equation}
\begin{equation}\label{Sij}
\begin{split}
\bar{\Lambda}_{ij}^{\pi,W}(n) & = \frac{\Lambda_{i,W}^\pi(n)}{\sum_{k=-1}^{n-1}\pi_k \sup_{\theta_j\in\Theta_j} LR_{j,\theta_j}(k,n)} ,
\\
& \quad i,j \in \Nc, ~ i \neq j, ~~ n \ge 1;
\\
\bar{\Lambda}_{i0}^{\pi,W}(n) & = \frac{\Lambda_{i,W}^\pi(n)}{\Pb(\nu \ge n)}, \quad i \in \Nc, ~~ n \ge 1 ,
\end{split}
\end{equation}
where in the statistic $\Lambda_{i, W}^\pi(n)$ defined in \eqref{LRmixturespi} $i=0$ corresponds to the hypothesis $\Hyp_0$ that there is no change (in the first $n$ observations).

Write $\Nc_0= \{0,1,\dots,N\}$. For the set of positive thresholds $A=(A_{ij})$, $j \in \Nc_0\setminus\{i\}$, $i \in \Nc$, 
the change detection--identification rule $\delta_A=(d_A, T_A)$ is defined as follows:
\begin{equation}
T_A = \min_{\ell\in \Nc} T_A^{(\ell)},  ~\quad d_A = i ~~ \text{if}~~ T_A= T_A^{(i)} , \label{CPDisolRule}
\end{equation}
where the Markov times $T_A^{(i)}$, $i \in \Nc$ are given by
\begin{equation}  \label{TAi}
T_A^{(i)} = \inf\set{n \ge 1: \bar{\Lambda}_{i j}^{\pi,W}(n) \ge A_{ij} ~ \forall~ j \in \Nc_0\setminus \{i\}} .
\end{equation}
 In definitions of stopping times we always set $\inf\{\varnothing\}=\infty$, i.e., $T_A^{(i)}=\infty$ if there is no such $n$.
If $T_A= T_A^{(i)}$ for several values of $i$ then any of them can be taken.

%__________________________________________________________________
\section{Optimization Problems and Assumptions} \label{sec:Ch8Problems}

Let $\Eb_{k,i,\theta_i}$ and $\Eb_\infty$ denote expectations under probability measures $\Pb_{k,i,\theta_i}$ and $\Pb_\infty$, respectively, where $\Pb_{k,i,\theta_i}$ corresponds to 
model \eqref{modeldetisol} with an assumed value of the parameter $\theta_i\in\Theta_i$, change point $\nu=k$, and the affected stream $i\in\Nc$. Define the probability measure 
$\Pb^\pi_{i,\theta_i} (\Ac\times \mc{K})=\sum_{k\in \mc{K}}\,\pi_{k} \Pb_{k,i,\theta_i}\left(\Ac\right)$  under which the change point $\nu$ has distribution $\pi=\{\pi_k\}_{k\in\Zbb}$ and the model for the observations
is of the form \eqref{modeldetisol} and let $\Eb^\pi_{i,\theta_i}$ denote the corresponding expectation.

For $r \ge 1$, $\nu=k \in \Zbb$, $\theta_i\in\Theta_i$, and $i\in\Nc$  introduce the risk associated with the conditional $r$th moment of the detection delay 
\begin{equation} \label{Riskdefki}
\Rc^r_{k,i,\theta_i}(\delta)=  \Eb_{k, i,\theta_i}\left[(T-k)^r;  d=i  \,|\, T> k \right], 
\end{equation}
where for $k=-1$ we set $T-k=T$, but not $T+1$, as well as the integrated (over prior $\pi$) risk associated with the moments of 
delay to detection  
\begin{equation} \label{Ch8Riskdefmultiple}
\begin{split}
\Rca^r_{i,\theta_i}(\delta) & : = \Eb^\pi_{i,\theta_i} [ (T-\nu)^r ;  d=i  \,|\, T> \nu]
\\
& = \frac{\Eb^\pi_{i,\theta_i} [ (T-\nu)^r,  d=i, T> \nu]}{\Pb_{i,\theta_i}^\pi(T>\nu)}
\\
& = \frac{{\displaystyle\sum_{k=-1}^\infty} \pi_k \Eb_{k, i,\theta_i}\left[(T-k)^r,  d=i, T> k \right]}{\sum_{k=-1}^\infty \pi_k \Pb_{k,i,\theta_i}(T>k)} 
\\
&=\frac{{\displaystyle\sum_{k=-1}^\infty} \pi_k \Rc^r_{k,i,\theta_i}(\delta) \Pb_{k,i,\theta_i}(T>k)}{\Pb(\nu \le 0)+ \sum_{k=1}^\infty \pi_k \Pb_{\infty}(T>k)} 
\\
&= 
\frac{{\displaystyle\sum_{k=-1}^\infty} \pi_k \Rc^r_{k,i,\theta_i}(\delta)\Pb_\infty( T >k)}{1-\PFA^\pi( \delta)} ,
\end{split}
\end{equation}
where 
\begin{equation} \label{PFAdef}
\begin{split}
\PFA^\pi(\delta) & =\Pb^\pi_{i,\theta_i}( T \le \nu) 
\\
& = \sum_{k=-1}^\infty \pi_k \, \Pb_{k,i,\theta_i} (T\le k)
\\
& = \sum_{k=0}^\infty \pi_k \, \Pb_\infty( T \le k) 
\end{split}
\end{equation}
is the weighted probability of false alarm.
Note that in \eqref{Ch8Riskdefmultiple} and \eqref{PFAdef} we used the equality $\Pb_{k,i,\theta_i}(T\le k)=\Pb_{\infty}(T\le k)$ since the event $\{T \le k\}$ belongs to the sigma-algebra 
$\Fc_k=\sigma(\Xb^k)$ and, hence, depends only on the first $k$ observations which distribution corresponds to the measure $\Pb_\infty$. This implies, in particular, that
\begin{align*}
\Pb_{i,\theta_i}^\pi(T>\nu) & = 1- \PFA^\pi(\delta)
\\
& = \Pb(\nu \le 0)+ \sum_{k=1}^\infty \pi_k \Pb_{\infty}(T>k).
\end{align*}

Also, introduce 
\begin{equation} \label{PFAidef}
\begin{split}
\PFA_i^\pi(\delta) & =\Pb^\pi_{i,\theta_i}( T \le \nu; d=i)
\\
& = \sum_{k=0}^\infty \pi_k \, \Pb_{i,\theta_i}( T \le k; d=i)
\\
& = \sum_{k=0}^\infty \pi_k \, \Pb_{\infty}( T \le k; d=i) ,
\end{split}
\end{equation}
the weighted probability of false alarm on the event $\{d=i\}$, i.e., the probability of raising the alarm with the decision $d=i$ that there is a change in the $i$th stream when there is no change.  

The loss associated with wrong identification is reasonable to measure by the maximal probabilities of wrong decisions (misidentification)
\begin{equation}\label{PMIsupdef}
\PMI_{ij}^\pi(\delta)= \sup_{\theta_i\in\Theta_i} \Pb_{i, \theta_i}^\pi(d=j; T< \infty| T>\nu), 
\end{equation}
$ i,j = 1,\dots,N, ~~ i \neq j$. Note that
\[
\begin{split}
&\Pb_{i, \theta_i}^\pi(d=j; T< \infty| T>\nu) =  \frac{\Pb_{i, \theta_i}^\pi(d=j; \nu < T< \infty)}{\Pb_{i,\theta_i}^\pi(T>\nu)}
\\
&= \frac{\sum_{k=-1}^\infty \pi_k \, \Pb_{k, i, \theta_i}(d=j; k < T< \infty)}{1-\PFA^\pi(\delta)}.
\end{split}
\]

Define the class of change detection--identification rules $\delta$ with constraints on the probabilities of false alarm $\PFA_i^\pi(\delta)$ and the probabilities of misidentification $\PMI_{ij}^\pi(\delta)$:
\begin{equation}\label{classdetisol}
\begin{split}
& \class_\pi(\alphab, \betab) =  \{\delta: \PFA_i^\pi(\delta) \le \alpha_i, ~ i\in\Nc ~~ \text{and}
\\
& \quad  \PMI_{ij}^\pi(\delta) \le \beta_{ij}, ~i,j  \in \Nc, i \neq j\},
\end{split}
\end{equation}
where $\alphab=(\alpha_1,\dots,\alpha_N)$ and $\betab=(\beta_{ij})_{i,j\in\Nc, i\neq j}$ are the sets of prescribed probabilities  $\alpha_i\in(0,1)$ and $\beta_{ij}\in(0,1)$.

Ideally, we would be interested in finding an optimal rule $\delta_{\rm opt}=(d_{\rm opt}, T_{\rm opt})$ that solves the optimization problem
\begin{equation*}
\Rca_{i,\theta_i}^r(\delta_{\rm opt}) = \inf_{\delta \in\class_\pi(\alphab, \betab)}\,\Rca_{i,\theta_i}^r(\delta) ~ \forall ~ \theta_i\in\Theta_i, ~ i \in \Nc.
\end{equation*}
However, this problem is intractable for arbitrary values of $\alpha_i\in(0,1)$ and $\beta_{ij}\in(0,1)$. For this reason, we will 
focus on the asymptotic problem assuming that the given probabilities $\alpha_i$ and $\beta_{ij}$ approach zero. 
To be more specific, we will be interested in proving that the proposed detection--identification rule $\delta_A=(d_A, T_A)$ defined in \eqref{CPDisolRule}--\eqref{TAi}
is first-order uniformly asymptotically optimal in the following sense
\begin{equation}\label{FOAOdefdetisol}
\begin{split}
& \lim_{\alpha_{\max}\to0, \beta_{\max}\to0} \frac{\displaystyle\inf_{\delta\in\class_\pi(\alphab, \betab)}\Rca_{i,\theta_i}^r(\delta)}{\Rca_{i,\theta_i}^r(\delta_A)} =1   
\\
&\qquad \text{for all} ~ \theta_i\in\Theta_i ~ \text{and} ~ i\in \Nc,
\end{split}
\end{equation} 
where $A=A(\alphab,\betab)$ is the set of suitably selected thresholds such that $\delta_A\in \class_\pi(\alphab, \betab)$. Hereafter $\alpha_{\max}=\max_{i\in\Nc} \alpha_i$, $\beta_{\max}= \max_{i,j\in\Nc, i\neq j} \beta_{ij}$.

In addition, we will prove that the rule $\delta_A=(d_A, T_A)$ is uniformly pointwise first-order asymptotically optimal in a sense of minimizing the conditional risk \eqref{Riskdefki} for all change point values 
$\nu=k\in \Zbb$, i.e.,
\begin{equation}\label{FOAOunifdef}
\begin{split}
& \lim_{\alpha_{\max}\to0, \beta_{\max}\to0} \frac{\displaystyle\inf_{\delta\in\class_\pi(\alphab, \betab)}\Rc_{k, i, \theta_i}^r(\delta)}{\Rc_{k, i, \theta_i}^r(\delta_A)} =1   
\\
& \quad \text{for all} ~ \theta_i\in \Theta_i, ~ k\in \Zbb, ~  i\in \Nc.
\end{split}
\end{equation} 

It is also of interest to consider the class of detection--identification rules 
\begin{equation}\label{classdetisol2}
\class_\pi^\star(\alpha, \bar{\betab})= \set{\delta: \PFA^\pi(\delta) \le \alpha, ~ \PMI_i^\pi(\delta) \le \bar{\beta}_i, ~ i\in \Nc}
\end{equation}
 ($\bar{\betab} =(\bar{\beta}_1,\dots,\bar{\beta}_N)$) with constrains on the total probability of false alarm  $\PFA^\pi(\delta)$ (defined in \eqref{PFAdef})
regardless of the decision $d=i$ which is made under hypothesis $\Hyp_\infty$ and on the misidentification probabilities
\[
\PMI_{i}^\pi(\delta)= \sup_{\theta_i\in\Theta_i} \Pb_{i, \theta_i}^\pi(d \neq i; T< \infty| T>\nu), \quad i\in \Nc.
\]
Obviously, $\PFA^\pi(\delta)= \sum_{i=1}^N \PFA_i^\pi(\delta)$ and $\PMI_{i}^\pi(\delta)= \sum_{j\in \Nc\setminus\{i\}}\PMI_{i j}^\pi(\delta)$.

In this paper, we consider only a fixed number of hypotheses $N$. The large-scale (Big Data) case where $N\to\infty$ with a certain rate (which requires a different definition of false alarm and 
misidentification rates) will be considered elsewhere.

 In the following, we assume that mixing measures $W_i$, $i=1,\dots,N$,  satisfy the condition:
\begin{equation}\label{CondonWi}
\begin{split}
&W_i\{\vartheta\in\Theta_i\,:\,\vert \vartheta-\theta_i\vert<\varkappa\}>0 
\\
&\qquad \text{for any}~ \varkappa>0 ~ \text{and any}~ \theta_i\in \Theta_i.
\end{split}
\end{equation}

By \eqref{LRdeteisol}, for the assumed values of $\nu = k$,  $i \in \Nc$, and $\theta_i\in\Theta_i$, the log-likelihood ratio (LLR) $\lambda_{i,\theta_i}(k, k+n)=\log LR_{i,\theta_i}(k,k+n)$ 
of observations accumulated by the time $k+n$ is 
\[
\lambda_{i,\theta_i}(k, k+n)= \sum_{t=k+1}^{k+n}\log \Lc_{i,\theta_i}(t), \quad n \ge 1,
\]
and the LLR between the hypotheses $\Hyp_{k,i}$ and $\Hyp_{k,j}$ of observations accumulated by the time $k+n$ is
\[
\begin{split}
& \lambda_{i,\theta_i; j, \theta_j}(k, k+n) = \log \frac{p(\Xb^{k+n}|\Hyp_{k, i})}{p(\Xb^{k+n}|\Hyp_{k, j})}  
\\
& \equiv \lambda_{i,\theta_i}(k, k+n) - \lambda_{j,\theta_j}(k, k+n), \quad n \ge 1.
\end{split}
\]
For $j=0$,  we set $\lambda_{0,\theta_0}(k, k+n)=0$, so that $\lambda_{i,\theta_i; 0, \theta_0}(k, k+n)= \lambda_{i,\theta_i}(k, k+n)$.

To study asymptotic optimality we need certain constraints imposed on the prior distribution $\pi=\{\pi_k\}$ and on the asymptotic behavior of the decision statistics as the sample size increases (i.e., on the general
stochastic model). 

For $\varkappa>0$, let $\Gamma_{\varkappa,\theta_i}=\{\vartheta\in\Theta_i\,:\,\vert \vartheta-\theta_i\vert<\varkappa\}$ and for $0<I_{i j}(\theta_i,\theta_j)<\infty$, $j\in\Nc_0\setminus\{i\}$,
$i\in\Nc$, define
\begin{align}
& p_{M,k}(\varepsilon; i,\theta_i; j, \theta_j)  = \nonumber
\\
&  \Pb_{k,i,\theta_i}\set{\frac{1}{M}\max_{1 \le n \le M} \lambda_{i,\theta_i; j, \theta_j}(k, k+n) \ge (1+\varepsilon) I_{i j}(\theta_i, \theta_j)}, \nonumber
\\
& \Upsilon_r(\varkappa, \varepsilon; i, \theta_i) = \nonumber
\\
&\sum_{n=1}^\infty \, n^{r-1} \, \sup_{k \in \Zbb_+} \Pb_{k,i,\theta_i}\Big\{\frac{1}{n} \inf_{\vartheta\in \Gamma_{\varkappa,\theta_i}}\lambda_{i,\vartheta}(k, k+n) \nonumber
\\
&< I_{i}(\theta_i)  - \varepsilon\Big\}  ,
\label{UpsilonCh8}
\end{align}
where $I_{i 0}(\theta_i, \theta_0)=I_i(\theta_i)$, so that
\[
\begin{split}
& p_{M,k}(\varepsilon; i,\theta_i; 0, \theta_0)  =p_{M,k}(\varepsilon; i,\theta_i)   
\\
& = \Pb_{k,i,\theta_i}\set{\frac{1}{M}\max_{1 \le n \le M} \lambda_{i,\theta_i}(k, k+n) \ge (1+\varepsilon) I_{i}(\theta_i)}.
\end{split}
\]

Regarding the model for the observations \eqref{modeldetisol}, we assume that the following two conditions are satisfied (for local LLRs in data streams): 
\vspace{2mm}

\noindent $\C_{1}$. {\em  There exist positive and finite numbers $I_{i}(\theta_i)$, $\theta_i\in \Theta_i$, $i \in \Nc$ and 
$I_{ij}(\theta_i,\theta_j)$, $\theta_j \in \Theta_j$, $ j \in \Nc\setminus\{i\}$, $\theta_i\in \Theta_i$, $i \in \Nc$, such that for any  $\varepsilon >0$}
\begin{equation}\label{Ch8Pmaxi}
\begin{split}
&\lim_{M\to\infty} p_{M,k}(\varepsilon; i,\theta_i; j, \theta_j) =0 
\\
& \text{for all}~ k\in \Zbb_+, ~  \theta_i\in\Theta_i, ~ \theta_j \in \Theta_j, ~ j \in \Nc_0\setminus\{i\}, ~ i \in \Nc  .
\end{split}
\end{equation}

\noindent $\C_{2}$. {\em  For any $\varepsilon>0$ and some $r\ge 1$}
\begin{equation}\label{Ch8rcompLefti}
\begin{split}
\lim_{\varkappa\to0}\Upsilon_r(\varkappa, \varepsilon; i,\theta_i) & < \infty \quad \text{for all}~  \theta_i\in\Theta_i, ~ i\in \Nc  .
\end{split}
\end{equation}

Note that condition $\C_{1}$  holds whenever $ \lambda_{i,\theta_i; j, \theta_j}(k, k+n)/n$ converges almost surely (a.s.) to $I_{ij}(\theta_i,\theta_j)$ under $\Pb_{k, i,\theta_i}$, i.e., 
for all $\theta_i\in \Theta_i$
\begin{equation}\label{sec:MaRe.1}
\Pb_{k, i, \theta_i}\set{\lim_{n\to\infty} \frac{1}{n}\lambda_{i,\theta_i; j, \theta_j}(k, k+n) = I_{ij}(\theta_i,\theta_j)} =1 .
\end{equation}

Regarding the prior distribution $\pi_k=\Pb(\nu=k)$ we assume that it is fully supported (i.e., $\pi_k>0$ for all $k\in\Zbb_+$, $0 \le \pi_{-1} <1$ and $\pi_\infty=0$) and 
the following two conditions are satisfied:
\vspace{2mm}

\noindent $\CP_1$. {\em  For some} $0 \le \mu <\infty$,
\begin{equation}\label{Prior}
\lim_{n\to\infty}\frac{1}{n}\abs{\log \sum_{k=n+1}^\infty \pi_k} = \mu.
\end{equation}
\noindent $\CP_2$. {\em If $\mu=0$, then in addition}
\begin{equation}\label{Prior1}
\sum_{k=0}^\infty \pi_k |\log\pi_k|^r < \infty \quad  \text{for some} ~ r\ge 1.
\end{equation}
The class of prior distributions satisfying conditions $\CP_1$ and $\CP_2$  will be denoted by $\Cb(\mu)$. 

Note that if $\mu >0$, then the prior distribution has an exponential right tail.  In this case, condition \eqref{Prior1} holds automatically. 
If $\mu=0$, the distribution has a heavy tail, i.e., belongs to the model with a vanishing hazard rate.
However, we cannot allow this distribution to have a too heavy tail, which
is guaranteed by condition  $\CP_2$. 
A typical heavy-tailed prior distribution that satisfies both conditions $\CP_1$ with $\mu=0$ and  $\CP_2$ for all $r\ge 1$ is a discrete Weibull-type distribution with the 
shape parameter $0<\kappa < 1$. Constraint \eqref{Prior1} is often guaranteed by finiteness of the $r$-th moment, $\Eb[\nu^r]<\infty$. 

 To obtain lower bounds for moments of the detection delay we need only right-tail conditions \eqref{Ch8Pmaxi}. However, to establish the asymptotic optimality property of the rule $\delta_A$ both 
right-tail and left-tail conditions \eqref{Ch8Pmaxi} and \eqref{Ch8rcompLefti} are needed.

%%_____________________________________________________________________________________
\section{Upper Bounds on Probabilities of False Alarm and Misidentification of the Detection--Identification Rule $\delta_A$}\label{sec:PFAandPMI}

Let $\wtPb_{i,\theta_i}^{\pi,n}(\Ac)=\Pb_{i,\theta_i}^\pi(\Ac, \nu<n)$ denote the measure $\Pb_{i,\theta_i}^\pi$ on the event $\{\nu < n\}$. 
The distribution $\wtPb_{i,\theta_i}^{\pi,n}(\Xb^n\in \Xc^n)$ has density
\begin{equation*}
\begin{split}
 f_{i,\theta_i}^{\pi,n}(\Xb^n)  & =  \sum_{k=-1}^{n-1} \Bigg [ \pi_k \prod_{t=1}^k g_i(X_t(i)|\Xb^{t-1}(i)) 
 \\
 & \prod_{t=k+1}^{n}  f_{i,\theta_i}(X_t(i)|\Xb^{t-1}(i)) \Bigg ]\times
\\
& \quad \Pb(\nu < n) \prod_{j\in \Nc\setminus \{i\}} \prod_{t=1}^n g_j(X_t(j)|\Xb^{t-1}(j)) ,
\end{split}
\end{equation*}
where $\prod_{t=1}^{-1} g_i(X_t(i)|\Xb^{t-1}(i))=1$. Write
\[
 f_{i,W}^{\pi,n}(\Xb^n) = \int_{\Theta_i}   f_{i,\theta_i}^{\pi,n}(\Xb^n) \, \mrm{d} W_i(\theta_i).
\]

Next, define the statistic $\widetilde{\Lambda}_{i,j,\theta_j}^{\pi,W}(n) = \Lambda_{i,W}^\pi(n)/\Lambda_{j,\theta_j}^\pi(n)$ and the measure 
\[
\wtPb_{\ell,W}^{\pi,n}(\Ac) = \int_{\Theta_\ell} \wtPb_{\ell,\theta_\ell}^{\pi,n}(\Ac) \mrm{d} W_\ell(\theta_\ell).
\]
Denote by $\Pb |_{\Fc_n}$ the restriction of the measure $\Pb$ to the sigma-algebra $\Fc_n$. Obviously,
\[
\widetilde{\Lambda}_{i,j,\theta_j}^{\pi,W}(n) = \frac{\mrm{d} \wtPb_{i,W}^{\pi,n} }{\mrm{d} \wtPb_{j,\theta_j}^{\pi,n}} \Bigg |_{\Fc_n}, \quad i \neq j,
\]
and hence, the statistic $\widetilde{\Lambda}_{i,j,\theta_j}^{\pi,W}(n)$ is a $(\wtPb_{j,\theta_j}^{\pi,n}, \Fc_n)$-martingale with unit expectation for all $\theta_j\in\Theta_j$. Therefore, by the Wald--Doob identity, 
for any stopping time $T$ and all $ \theta_j\in\Theta_j,$
\begin{equation}\label{WDidentity}
\begin{split}
& \wtE_{i,\theta_i}^\pi \brcs{\widetilde{\Lambda}_{j,i,\theta_i}^{\pi,W}(T)\Ind{\Ac, T<\infty}}  = \wtE_{j, W}^\pi \brcs{\Ind{\Ac, T<\infty}} 
\\
& = \wtPb_{j,W}^{\pi,T}(\Ac\cap \{T<\infty\}),
\end{split}
\end{equation}
where $\wtE_{j, W}^\pi$ and  $\wtE_{j, \theta_j}^\pi$  stand for the operators of expectation under $\wtPb_{j,W}^{\pi,T}$ and $\wtPb_{j,\theta_j}^{\pi,T}$, respectively.

The following theorem establishes upper bounds for the PFA and PMI of the proposed detection--identification rule $\delta_A$. Note that these bounds are valid in the most general case -- neither of the conditions 
on the model $\Cb_1$, $\Cb_2$ or on the prior distribution $\CP_1$, $\CP_2$ are required.

%%Theorem
\begin{theorem}\label{Th:UpperboundsPFAPMI}
Let $\delta_A$ be the changepoint detection--identification rule defined in  \eqref{CPDisolRule}--\eqref{TAi}. The following upper bounds for the PFA and PMI of rule $\delta_A$ hold
\begin{align}
\PFA_i^\pi(\delta_A) & \le (1+A_{i0})^{-1}, \quad i \in \Nc, \label{UpperPFAi}
\\
\PFA^\pi(\delta_A) & \le \sum_{i=1}^N (1+A_{i0})^{-1} \label{UpperPFAdetisol}
\end{align}
and
\begin{align}
\PMI_{ij}^\pi(\delta_A) & \le \frac{1+A_{i0}}{A_{i0} \, A_{ji}}, \quad i, j \in \Nc, ~ i \neq j , \label{UpperPMIij}
\\
\PMI_{i}^\pi(\delta_A) & \le \frac{1+A_{i0}}{A_{i0}} \sum_{j \in \Nc\setminus\{i\}}\frac{1}{A_{ji}}, \quad i \in \Nc . \label{UpperPMIi}
\end{align}
Thus, if $\alpha_{\max} < 1-\pi_{-1}$, then
\begin{equation}\label{Aalpha}
\begin{split}
& A_{i0}  = \frac{1-\alpha_i}{\alpha_i} ~~ \text{and}~~ A_{ij} =\frac{1}{(1-\alpha_j)\beta_{ji}}  
\\
& \quad\text{imply} ~~ \delta_A\in \class_\pi(\alphab,\betab) ,
\end{split}
\end{equation}
and if $A_{i0}=A_0$ for $i\in \Nc$ and $A_{ij}=A_j$ for $j \in \Nc\setminus\{i\}$, then
\begin{equation}\label{Aalpha1}
\begin{split}
& A_{0} = \frac{N}{\alpha} (1-\alpha/N) ~~ \text{and}~~ A_{j} =\frac{N-1}{(1-\alpha/N)\bar{\beta}_j}  
\\
& \quad\text{imply} ~~ \delta_A\in \class_\pi^\star(\alpha, \bar{\betab}).
\end{split}
\end{equation}
\end{theorem}

%%Proof
\begin{IEEEproof}
Using the Bayes rule, notation \eqref{LRdeteisol}--\eqref{Sij}, and the fact that $LR_{i,\theta_i}(k,n) =1$ for $k \ge n$, we obtain   
\begin{align*}
\Pb(\nu=k|\Fc_n)  & =   \frac{\pi_k LR_{i,W}(k,n)}{\sum_{j=-1}^\infty \pi_j LR_{i,W}(j,n)}
\\
= &  \frac{\pi_k LR_{i,W}(k,n)}{\sum_{j=-1}^{n-1} \pi_j LR_{i,W}(j,n) + \Pb(\nu\ge n)},
\end{align*}
so that 
\begin{align*}
\Pb(\nu \ge n | \Fc_n) & = \sum_{k=n}^\infty  \Pb(\nu=k|\Fc_n) 
 =  \frac{\Pb(\nu\ge n)}{\Lambda_{i,W}^\pi + \Pb(\nu\ge n)}
 \\
 & = \frac{1}{\bar{\Lambda}_{i0}^{\pi,W}(n) + 1} .
\end{align*}
Next, obviously,
\[
\PFA_i(\delta_A) = \Pb_{i,\theta_i}^\pi(T_A^{(i)} \le \nu, T_A=T_A^{(i)}) \le \Pb_{i,\theta_i}^\pi(T_A^{(i)} \le \nu).
\]
Therefore, taking into account that $\Pb_{i,\theta_i}^\pi(T_A^{(i)} \le \nu) = \Eb_{i,\theta_i}^\pi [\Pb(T_A^{(i)} \le \nu | \Fc_{T_A^{(i)}})]$
and that $\bar{\Lambda}_{i0}^{\pi,W}(T_A^{(i)}) \ge A_{i0}$ on $\{T_A^{(i)}<\infty\}$, we obtain
\[
\begin{split}
\PFA_i(\delta_A)  &\le \Eb_{i,\theta_i}^\pi [(1+ \bar{\Lambda}_{i0}^{\pi,W}(T_A^{(i)}))^{-1}; T_A^{(i)}< \infty] 
\\
&\le 1/(1+A_{i0})
\end{split}
\]
and inequalities \eqref{UpperPFAi} follow. Inequality \eqref{UpperPFAdetisol} follows immediately from the fact that $\PFA^\pi(\delta)= \sum_{i=1}^N \PFA_i^\pi(\delta)$. 

To prove the upper bound \eqref{UpperPMIij} note that $\widetilde{\Lambda}_{j, i, \theta_i}^{\pi,W}(n) \ge \bar{\Lambda}_{ji}^{\pi,W}(n)$ for all $n \ge 1$ and $\theta_i\in\Theta_i$ and that 
$\bar{\Lambda}_{ji}^{\pi,W}(T_A^{(j)}) \ge A_{ji}$ on $\{T_A^{(j)}< \infty\}$ and we have
\begin{align*}
& \Pb_{i,\theta_i}^\pi(d_A=j, \nu < T_A <\infty)  
\\
&= \Pb_{i,\theta_i}^\pi(T_A= T_A^{(j)}, \nu < T_A^{(j)} <\infty) 
\\
&= \wtPb_{i,\theta_i}^{\pi,T_A}(T_A= T_A^{(j)}, T_A^{(j)} <\infty)
\\
& =\wtE_{i,\theta_i}^\pi\brcs{\frac{\bar{\Lambda}_{ji}^{\pi,W}(T_A^{(j)})}{\bar{\Lambda}_{ji}^{\pi,W}(T_A^{(j)})}\Ind{T_A= T_A^{(j)},T_A^{(j)}<\infty}}
\\
&\le \frac{1}{A_{ji}}\wtE_{i,\theta_i}^\pi\brcs{\bar{\Lambda}_{ji}^{\pi,W}(T_A^{(j)})\Ind{T_A= T_A^{(j)},T_A^{(j)} <\infty}}
\\ 
& \le \frac{1}{A_{ji}}\wtE_{i,\theta_i}^\pi\brcs{\widetilde{\Lambda}_{j, i, \theta_i}^{\pi,W}(T_A^{(j)})\Ind{T_A= T_A^{(j)},T_A^{(j)} <\infty}} ~ \forall~ \theta_i\in\Theta_i,
\end{align*}
where, by equality \eqref{WDidentity}, the last term is equal to
\[
\frac{1}{A_{ji}}\wtPb_{j,W}^{\pi,T_A}(\{T_A= T_A^{(j)}\}\cap \{T_A^{(j)}<\infty\}). 
\]
This yields
\begin{align*}
& \Pb_{i,\theta_i}^\pi(d_A=j, \nu < T_A <\infty) 
 \\
 &\le \frac{1}{A_{ji}}\wtPb_{j,W}^{\pi,T_A}(\{T_A= T_A^{(j)}\}\cap \{T_A^{(j)}<\infty\}) 
 \\
 &\le \frac{1}{A_{ji}} ~~ \text{for all}~ \theta_i\in\Theta_i.
\end{align*}
Since $\Pb_{i,\theta_i}^\pi(d_A=j | T_A> \nu) =\Pb_{i,\theta_i}^\pi(d_A=j, \nu < T_A <\infty)/\Pb_{i,\theta_i}^\pi(T_A>\nu)$ and, by \eqref{UpperPFAi}, $\Pb_{i,\theta_i}^\pi(T_A>\nu) \ge  A_i/(1+A_i)$, the upper bound
\eqref{UpperPMIij} follows. The upper bound \eqref{UpperPMIi} follows from \eqref{UpperPMIij} and the fact that $\PMI_{i}^\pi(\delta)= \sum_{j\in\Nc\setminus\{i\}}\PMI_{i j}^\pi(\delta)$.

Implications \eqref{Aalpha} and  \eqref{Aalpha1} are obvious.
\end{IEEEproof}

%%Remark
\begin{remark}\label{RemUBPFA}
Typically, the upper bounds \eqref{UpperPFAi}--\eqref{UpperPMIi} for PFA and PMI are not tight but rather quite conservative, especially when overshoots over thresholds are large (i.e., when the hypotheses 
$\Hyp_i$ and $\Hyp_\infty$ are not close). Unfortunately, in the general non-i.i.d.\ case, the improvement of these bounds is not possible. In the i.i.d.\ case where observations are independent and identically distributed
with the common pre-change density $g_i(x)$ and the common post-change density $f_i(x)$ in the $i$th stream (i.e., when the post-change hypotheses are simple), it is possible to obtain asymptotically 
accurate approximations using the renewal theory similarly to how it was done in \cite[Th 7.1.5, p. 327]{TNB_book2014} for the PFA in the single-stream case.
\end{remark}

%%_________________________________________________________-
\section{Lower Bounds on the Moments of the Detection Delay in Classes $\class_\pi(\alphab,\betab)$ and $\class_\pi^\star(\alpha,\bar\betab)$}\label{sec:LBdetisol}

For $i\in \Nc$, define
\begin{equation}\label{Psii}
\begin{split}
\Psi_i(\alphab, \betab)  = \max\Bigg\{\frac{|\log \alpha_i|}{I_i(\theta_i)+\mu}, 
 \max_{j\in\Nc\setminus\{i\}} \frac{|\log\beta_{ji}|}{\displaystyle{\inf_{\theta_j\in\Theta_j} I_{ij}(\theta_i, \theta_j)}}\Bigg\}
\end{split}
\end{equation}
and
\begin{equation}\label{Psii2}
\begin{split}
\Psi_i^\star(\alpha, \bar{\betab}) = \max\Bigg\{\frac{|\log \alpha|}{I_i(\theta_i)+\mu}, 
\max_{j\in\Nc\setminus\{i\}} \frac{|\log\bar{\beta}_{j}|}{\displaystyle{\inf_{\theta_j\in\Theta_j} I_{ij}(\theta_i, \theta_j)}}\Bigg\}.
\end{split}
\end{equation}

The following theorem establishes asymptotic lower bounds on moments of the detection delay  $\Rc^r_{k, i,\theta_i}(\delta)$ 
and $\Rca^r_{i,\theta_i}(\delta)$ ($r\ge 1$)
in classes of detection--identification rules $\class_\pi(\alphab,\betab)$ and $\class_\pi^\star(\alpha,\bar{\betab})$ defined in \eqref{classdetisol} and \eqref{classdetisol2}, respectively. 
These bounds will be used in the next section for proving asymptotic optimality of the detection--identification rule $\delta_A$ with suitable thresholds. 

%%Theorem
\begin{theorem}\label{Th:LBdetisol}
Let, for some $\mu\ge 0$, the prior distribution belong to class $\Cb(\mu)$. 
Assume that for some positive and finite numbers $I_i(\theta_i)$ ($\theta_i\in\Theta_i$, $i\in \Nc$) and $I_{ij}(\theta_i,\theta_j)$ ($\theta_i\in \Theta_i$, $\theta_j\in\Theta_j$, $i\in \Nc$, $j\in\Nc\setminus\{i\}$) 
condition $\C_1$ holds. If $\inf_{\theta_j \in \Theta_j} I_{ij}(\theta_i, \theta_j) > 0$ for all $j \in \Nc\setminus \{i\}$, then for all $r >0$,  $\theta_i\in\Theta_i$, and $i\in\Nc$,
\begin{equation}\label{LBinclassdetisolk}
\liminf_{\alpha_{\max}, \beta_{\max}\to0} \frac{{\displaystyle\inf_{\delta\in\class_\pi(\alphab,\betab)}}  \Rc^r_{k, i,\theta_i}(\delta)}{[\Psi_{i}(\alphab, \betab)]^r} \ge  1 ~~ \text{for all} ~ k \in\Zbb,
\end{equation}
\begin{equation}\label{LBinclassdetisol}
\liminf_{\alpha_{\max}, \beta_{\max}\to0} \frac{{\displaystyle\inf_{\delta\in\class_\pi(\alphab,\betab)}}  \Rca^r_{i,\theta_i}(\delta)}{[\Psi_i(\alphab, \betab)]^r} \ge  1 
\end{equation}
and
\begin{equation}\label{LBinclassdetisolk2}
\liminf_{\alpha_{\max}, \beta_{\max}\to0} \frac{{\displaystyle\inf_{\delta\in\class_\pi^\star(\alpha,\bar{\betab})}}  \Rc^r_{k,i,\theta_i}(\delta)}{[\Psi_i^\star(\alpha, \bar{\betab})]^r} \ge  1 ~~ \text{for all} ~ k\in\Zbb ,
\end{equation}
\begin{equation}\label{LBinclassdetisol2}
\liminf_{\alpha_{\max}, \beta_{\max}\to0} \frac{{\displaystyle\inf_{\delta\in\class_\pi^\star(\alpha,\bar{\betab})}}  \Rca^r_{i,\theta_i}(\delta)}{[\Psi_i^\star(\alpha, \bar{\betab})]^r} \ge  1 ,
\end{equation}
where $\Psi_i(\alphab, \betab)$ and $\Psi_i^\star(\alpha, \bar{\betab})$ are defined in \eqref{Psii} and \eqref{Psii2}, respectively.
\end{theorem}

%%Proof
\begin{IEEEproof}
We only provide the proof of asymptotic lower bounds \eqref{LBinclassdetisolk} and \eqref{LBinclassdetisol}. The proof of \eqref{LBinclassdetisolk2}  and \eqref{LBinclassdetisol2} is essentially similar. 

Notice that the proof can be split into two parts since if we show that, on one hand, for any rule $\delta\in\class_\pi(\alphab, \betab)$  as $\alpha_{\max}\to0, \beta_{\max}\to0$
\begin{equation}\label{LBinclassdetisolkPMI}
\Rc^r_{k,i,\theta_i}(\delta) \ge \max_{j\in\Nc \setminus\{i\}} \brcs{\frac{|\log\beta_{ji}|}{\displaystyle{\inf_{\theta_j\in \Theta_j} I_{ij}(\theta_i, \theta_j)}}}^r (1+o(1)) ~~ \forall~k \in\Zbb
\end{equation}
and
\begin{equation}\label{LBinclassdetisolPMI}
\Rca^r_{i,\theta_i}(\delta) \ge \max_{j\in\Nc \setminus\{i\}} \brcs{\frac{|\log\beta_{ji}|}{\displaystyle{\inf_{\theta_j\in \Theta_j} I_{ij}(\theta_i, \theta_j)}}}^r (1+o(1)) ,
\end{equation}
and on the other hand
\begin{equation}\label{LBinclassdetisolkPFA}
 \Rc^r_{k,i,\theta_i}(\delta) \ge  \brc{\frac{|\log \alpha_i|}{I_i(\theta_i)+\mu}}^r (1+o(1)) ~~ \forall~k \in\Zbb,
\end{equation}
and
\begin{equation}\label{LBinclassdetisolPFA}
 \Rca^r_{i,\theta_i}(\delta) \ge  \brc{\frac{|\log \alpha_i|}{I_i(\theta_i)+\mu}}^r (1+o(1)),
\end{equation}
where $o(1)\to0$, then, obviously, combining inequalities \eqref{LBinclassdetisolkPMI} and \eqref{LBinclassdetisolkPFA} yields \eqref{LBinclassdetisolk} and combining 
\eqref{LBinclassdetisolPMI} and \eqref{LBinclassdetisolPFA} yields \eqref{LBinclassdetisol}.

The detailed proof of inequalities \eqref{LBinclassdetisolkPMI}--\eqref{LBinclassdetisolPFA} is postponed to the Appendix. 
\end{IEEEproof}

%%_____________________________________________
\section{Asymptotic Optimality}\label{sec:Asoptdetisol}

The following proposition, whose proof is given in the Appendix,  establishes first-order asymptotic approximations to the moments of the detection delay of the detection--identification rule $\delta_A$ 
when thresholds $A_{ij}$ go to infinity regardless of the PFA and PMI constraints.  Write $A_{\min} = \min_{i \in \Nc, j \in \Nc_0\setminus\{i\}} A_{ij}$.

%%% Proposition
\begin{proposition}\label{Pr:FOapproxAdetisol}
 Let $r\ge 1$ and let the prior distribution of the change point belong to class $\Cb(\mu)$.  Assume that for some  $0< I_i(\theta_i)<\infty$, $\theta_i\in\Theta_i$, $i\in \Nc$ and $0<I_{ij}(\theta_i,\theta_j)<\infty$, 
 $\theta_i\in \Theta_i$, $\theta_j\in\Theta_j$, $i\in \Nc$, $j\in\Nc\setminus\{i\}$  right-tail and left-tail conditions
$\C_{1}$ and $\C_{2}$ are satisfied and that $\inf_{\theta_j \in \Theta_j} I_{ij}(\theta_i, \theta_j) > 0$ for all $j\in \Nc\setminus\{i\}$, $i \in \Nc$. Then, for all $0<m \le r$, $\theta_i\in\Theta_i$, and $i \in \Nc$ 
 as $A_{\min} \to \infty$
 \begin{equation} \label{MADDAsaprk}
\Rc^m_{k, i,\theta_i}(\delta_A)  \sim \brcs{\Psi_i(A, \theta_i, \mu) }^m ~~ \text{for all}~ k \in \Zbb 
\end{equation}
and
\begin{equation} \label{MADDAsapr}
\Rca^m_{i,\theta_i}(\delta_A)  \sim \brcs{\Psi_i(A, \theta_i, \mu) }^m  ,
\end{equation}
where
\begin{equation}\label{PsiACh8}
\begin{split}
\Psi_i(A, \theta_i, \mu) =  \max \Bigg\{\frac{\log A_{i0}}{I_i(\theta_i)+\mu}, 
\max_{j\in\Nc\setminus\{i\}} \frac{\log A_{ij}}{\displaystyle{\inf_{\theta_j\in\Theta_j} I_{ij}(\theta_i, \theta_j)}}\Bigg\}.
\end{split}
\end{equation}
\end{proposition}
Hereafter we use a standard notation $x_a\sim y_a$ as $a \to a_0$ if $\lim_{a\to a_0}(x_a/y_a)=1$.

In order to prove this proposition we need the following lemma, whose proof is given in the Appendix. For $i=1,\dots,N$,  define 
\begin{align*}
\lambda_{i,W}(k,k+n) & = \log LR_{i,W}(k,k+n),
\\
\lambda_i^\pi(n) & = \log \brcs{\sum_{k=-1}^{n-1}\pi_k \sup_{\theta_i\in\Theta_i} LR_{i,\theta_i}(k,n)},
\end{align*}
\begin{equation}\label{PsiA}
\begin{split}
&\widetilde{\Psi}_i(A, \pi_k, \theta_i,\mu, \varepsilon) = \max\Bigg\{\frac{\log (A_{i0}/\pi_k)}{I_i(\theta_i)+\mu-\varepsilon}, 
\\
& \qquad \max_{j\in\Nc\setminus\{i\}} \frac{\log (A_{ij}/\pi_k)}{\inf_{\theta_j\in\Theta_j} I_{ij}(\theta_i, \theta_j) -\varepsilon}\Bigg\} ,
\end{split}
\end{equation}
\begin{align*}
M_{i}(A)& =M_{i}(A, \pi_k, \theta_i, \mu, \varepsilon)=1+\left\lfloor \widetilde{\Psi}_i(A, \pi_k, \theta_i,\mu, \varepsilon) \right \rfloor , \nonumber
\end{align*}
where $\left\lfloor y  \right \rfloor$ is the greatest integer.

%%Lemma
\begin{lemma}\label{Lem:UpperEki}
Let $r\ge 1$ and let the prior distribution of the change point satisfy condition \eqref{Prior}. Then, for a sufficiently large $A_{\min}$, any $0<\varepsilon <J_{ij}(\theta_i,\mu)$ and all $k\in\Zbb$,
\begin{equation}\label{EkineqCh8}
\begin{split}
&  \Eb_{k,i,\theta_i}[\brc{T_A-k}^+]^r  \le  \brcs{1+\widetilde{\Psi}_i(A, \pi_k, \theta_i,\mu, \varepsilon)}^r 
 \\
& + r 2^{r-1} \sum_{n=M_{i}(A)}^\infty n^{r-1}\Pb_{k,i, \theta_i}\Bigg\{\frac{1}{n} \inf_{\vartheta\in\Gamma_{\varkappa,\theta_i}}\lambda_{i,\vartheta}(k, k+n) 
\\
& < I_{i}(\theta_i)  - \varepsilon\Bigg\} ,
\end{split}
\end{equation}
where $T_A-k=T_A$ for $k=-1$, $x^+=\max(0,x)$, and
\[
J_{ij}(\theta_i,\mu) = \min\{I_i(\theta_i) + \mu, \min_{j\in \Nc\setminus\{i\}} \inf_{\theta_j \in \Theta_j} I_{ij}(\theta_i,\theta_j)\} .
\]
\end{lemma}

Theorem~\ref{Th:UpperboundsPFAPMI}, Theorem~\ref{Th:LBdetisol} and Proposition~\ref{Pr:FOapproxAdetisol} allow us to conclude that the detection--identification rule $\delta_A$ is asymptotically 
first-order optimal in classes $\class_\pi(\alphab,\betab)$ and $\class_\pi^\star(\alpha, \bar{\betab})$ as $\alpha_{\max}, \beta_{\max} \to 0$.

%%Theorem
\begin{theorem}\label{Th:FOAOCh8} 
Let $r\ge 1$ and let the prior distribution of the change point belong to class $\Cb(\mu)$.  Assume that for some  $0< I_i(\theta_i)<\infty$, $\theta_i\in\Theta_i$, $i\in \Nc$ and $0<I_{ij}(\theta_i,\theta_j)<\infty$, 
$\theta_i\in \Theta_i$, $\theta_j\in\Theta_j$, $i\in \Nc$, $j\in\Nc\setminus\{i\}$  right-tail and left-tail conditions
$\C_{1}$ and $\C_{2}$ are satisfied and that $\inf_{\theta_j \in \Theta_j} I_{ij}(\theta_i, \theta_j) > 0$ for all $j\in \Nc\setminus\{i\}$, $i \in \Nc$. 

\vspace{2mm}
\noindent{\bf (i)}  If thresholds $A_{i0}$, $i\in\Nc$ and $A_{ij}$, $j \in \Nc\setminus\{i\}$, $i \in \Nc$ are so selected that $\PFA_i^\pi(\delta_A) \le \alpha_i$, $\PMI_{ij}(\delta_A) \le \beta_{ij}$ and $\log A_{i0} \sim |\log\alpha_i|$, 
$\log A_{ij} \sim |\log \beta_{ji}|$ as  $\alpha_{\max}, \beta_{\max} \to 0$, in particular as $A_{i0}=(1-\alpha_i)/\alpha_i$ and  $ A_{ij} =[(1-\alpha_j)\beta_{ji}]^{-1}$, then $\delta_{A}$ is first-order asymptotically optimal as 
$\alpha_{\max}, \beta_{\max} \to 0$ in class $\class_\pi(\alphab,\betab)$,  minimizing moments of the detection delay up to order $r$: for all $0<m \le r$, $\theta_i\in\Theta_i$, and  $i\in\Nc$
\begin{equation}\label{FOAOmomentsgenkCh8}  
\begin{split}
&  \inf_{\delta \in \class_\pi(\alphab,\betab)} \Rc^m_{k, i,\theta_i}(\delta)    \sim \Rc^m_{k, i,\theta_i}(\delta_A)  
\\
&  \sim \max\Bigg\{\frac{|\log \alpha_{i}|}{I_i(\theta_i)+\mu}, \max_{j\in\Nc\setminus\{i\}} \frac{|\log \beta_{ji}|}{\inf_{\theta_j\in\Theta_j} I_{ij}(\theta_i, \theta_j)}\Bigg\}^m
 \\
& \quad \text{as}~ \alpha_{\max}, \beta_{\max} \to 0~~ \text{for all}~ k\in\Zbb
 \end{split}
 \end{equation}
and
\begin{equation}\label{FOAOmomentsgenCh8} 
\begin{split}
&  \inf_{\delta \in \class_\pi(\alphab,\betab)} \Rca^m_{i,\theta_i}(\delta)   \sim  \Rca^m_{i,\theta_i}(\delta_A) 
\\
& \sim \max\Bigg\{\frac{|\log \alpha_{i}|}{I_i(\theta_i)+\mu}, \max_{j\in\Nc\setminus\{i\}} \frac{|\log \beta_{ji}|}{\inf_{\theta_j\in\Theta_j} I_{ij}(\theta_i, \theta_j)}\Bigg\}^m
 \\
& \qquad   \text{as}~ \alpha_{\max}, \beta_{\max} \to 0.
 \end{split}
 \end{equation}

\noindent{\bf (ii)}  If thresholds $A_{i0}=A_0$ and $A_{ij}=A_{j}$, $j \in \Nc\setminus\{i\}$, $i \in \Nc$ are so selected that $\PFA^\pi(\delta_A) \le \alpha$, $\PMI_i(\delta_A) \le \bar{\beta}_i$ and $\log A_{0} \sim |\log\alpha|$, 
$\log A_{j} \sim |\log \bar{\beta}_{j}|$ as  $\alpha, \bar{\beta}_{\max} \to 0$, in particular as $A_{0}=N(1-\alpha/N)/\alpha$ and $ A_{j} =(N-1)[(1-\alpha/N)\bar{\beta}_{j}]^{-1}$, then $\delta_{A}$ is first-order asymptotically 
optimal as $\alpha, \bar{\beta}_{\max} \to 0$ in class $\class_\pi^\star(\alpha,\bar{\betab})$,  minimizing moments of the detection delay up to order $r$: for all $0<m \le r$, $\theta_i\in\Theta_i$, and  $i\in\Nc$,
    \begin{equation}\label{FOAOmomentsgenkCh8_2}  
\begin{split}
& \inf_{\delta \in \class_\pi^\star(\alpha,\bar{\betab})} \Rc^m_{k, i,\theta_i}(\delta)    \sim \Rc^m_{k, i,\theta_i}(\delta_A)   
\\
& \sim  \max\Bigg\{\frac{|\log \alpha|}{I_i(\theta_i)+\mu}, \max_{j\in\Nc\setminus\{i\}} \frac{|\log \bar{\beta}_{j}|}{\inf_{\theta_j\in\Theta_j} I_{ij}(\theta_i, \theta_j)}\Bigg\}^m
 \\
 & \qquad \text{as}~ \alpha, \bar{\beta}_{\max} \to 0 ~~ \text{for all} ~ k\in\Zbb
 \end{split}
 \end{equation}
and
\begin{equation}\label{FOAOmomentsgenCh8_2}
\begin{split}
& \inf_{\delta \in \class_\pi^\star(\alpha,\bar{\betab})} \Rca^m_{i,\theta_i}(\delta)  \sim \Rca^m_{i,\theta_i}(\delta_A)  
\\
& \sim  \max\Bigg\{\frac{|\log \alpha|}{I_i(\theta_i)+\mu}, \max_{j\in\Nc\setminus\{i\}} \frac{|\log \bar{\beta}_{j}|}{\inf_{\theta_j\in\Theta_j} I_{ij}(\theta_i, \theta_j)}\Bigg\}^m
 \\
 & \qquad  \text{as}~ \alpha, \bar{\beta}_{\max} \to 0 .
 \end{split}
 \end{equation}
\end{theorem}

%%%Proof
\begin{IEEEproof}
Proof of (i). Setting $\log A_{i0} \sim |\log\alpha_i|$ and $\log A_{ij} \sim |\log \beta_{ji}|$ in \eqref{MADDAsaprk} yields as $\alpha_{\max},\beta_{\max}\to 0$ 
\begin{equation}\label{AsapproxCh8}
\begin{split}
&\Rc^m_{k,i,\theta_i}(\delta_A)  \sim \max\Bigg\{\frac{|\log \alpha_{i}|}{I_i(\theta_i)+\mu}, 
\\
& \qquad \max_{j\in\Nc\setminus\{i\}} \frac{|\log \beta_{ji}|}{\inf_{\theta_j\in\Theta_j} I_{ij}(\theta_i, \theta_j)}\Bigg\}^m, ~~ i \in \Nc.
\end{split}
\end{equation}
In particular, $\log A_{i0} \sim |\log\alpha_i|$ and $\log A_{ij} \sim |\log \beta_{ji}|$ if $A_{i0}=(1-\alpha_i)/\alpha_i$ and  $ A_{ij} =[(1-\alpha_j)\beta_{ji}]^{-1}$, and by Theorem~\ref{Th:UpperboundsPFAPMI}, 
$\PFA_i^\pi(\delta_A) \le \alpha_i$ and $\PMI_{ij}(\delta_A) \le \beta_{ij}$ with this choice of thresholds (see~\eqref{Aalpha}).
Comparing asymptotic approximations \eqref{AsapproxCh8} with the lower bounds \eqref{LBinclassdetisolk}  in Theorem~\ref{Th:LBdetisol} completes the proof of \eqref{FOAOmomentsgenkCh8}. The proof of
\eqref{FOAOmomentsgenCh8} is similar.

Proof of (ii).  Setting $\log A_{i0}=\log A_{0} \sim |\log\alpha|$ and $\log A_{ij} =\log A_{j} \sim |\log \bar{\beta}_{j}|$ in \eqref{MADDAsaprk} yields as $\alpha_{\max},\bar{\beta}_{\max}\to 0$ 
\begin{equation}\label{AsapproxCh8_2}
\begin{split}
& \Rc^m_{k,i,\theta_i}(\delta_A) \sim \max\Bigg\{\frac{|\log \alpha|}{I_i(\theta_i)+\mu}, 
\\
& \qquad \max_{j\in\Nc\setminus\{i\}} \frac{|\log \bar{\beta}_j|}{\inf_{\theta_j\in\Theta_j} I_{ij}(\theta_i, \theta_j)}\Bigg\}^m, ~~ i \in \Nc.
\end{split}
\end{equation}
In particular, $\log A_{0} \sim |\log\alpha|$ and $\log A_j \sim |\log \bar{\beta}_{j}|$ if $A_{0}=N(1-\alpha/N)/\alpha$ and $ A_{j} =(N-1)[(1-\alpha/N)\bar{\beta}_{j}]^{-1}$, and by Theorem~\ref{Th:UpperboundsPFAPMI}, 
$\PFA^\pi(\delta_A) \le \alpha$ and $\PMI_{i}(\delta_A) \le \bar{\beta}_{i}$ with this choice of thresholds (see~\eqref{Aalpha1}).
Comparing asymptotic approximations \eqref{AsapproxCh8_2} with the lower bounds \eqref{LBinclassdetisolk2}  in Theorem~\ref{Th:LBdetisol} completes the proof of \eqref{FOAOmomentsgenkCh8_2}. 
The proof of \eqref{FOAOmomentsgenCh8_2} is similar.
\end{IEEEproof}

%%Remark
\begin{remark}
If the prior distribution $\pi=\pi^{\alpha_{\max},\beta_{\max}}$ depends on the PFA $\alpha_{\max}$ and PMI $\beta_{\max}$ constraints and  $\mu_{\alpha_{\max},\beta_{\max}} \to0$ as 
$\alpha_{\max},\beta_{\max}\to0$, then a modification of the preceding argument
can be used to show that the assertions of Theorem~\ref{Th:FOAOCh8} hold with $\mu=0$.
\end{remark}

Note that conditions \eqref{Ch8Pmaxi} are satisfied if
\[
\frac{1}{n}\lambda_{i,\theta_i; j, \theta_j}(k, k+n) \xra[n\to\infty]{\Pb_{k,i,\theta_i}-\text{a.s.}} I_{i j}(\theta_i, \theta_j)
\]
(see Lemma~B.1 in \cite[p. 243]{AGT_book2020}). Assume also that for some positive and finite numbers $ I_{0, i}(\theta_i)$, $i\in \Nc$,
\[
- \frac{1}{n}\lambda_{i,\theta_i}(k, k+n) \xra[n\to\infty]{\Pb_{\infty}-\text{a.s.}} I_{0, i}(\theta_i).
\]
In particular, in the i.i.d.\ case, these conditions hold with 
\begin{equation*}%\label{KLCh8}
\begin{split}
I_{i j}(\theta_i, \theta_j) & \equiv \Kc_{ij}(\theta_i,\theta_j) = \int \brc{\log \frac{f_{i, \theta_i}(x)}{f_{j,\theta_j}(x)}} f_{i, \theta_i}(x) \D x,  
\\
I_{0, i}(\theta_i)& \equiv \Kc_{0,i}(\theta_i) = \int \brc{\log \frac{g_{i}(x)}{f_{i,\theta_i}(x)}} g_{i}(x) \D x
\end{split}
\end{equation*}
being the Kullback--Leibler information numbers. Then, $I_{i j}(\theta_i, \theta_j) = I_i(\theta_i) + I_{0,j}(\theta_j) \ge I_i(\theta_i)$. Therefore, if the prior distribution of the change point is heavy-tailed (i.e., $\mu=0$) 
and the PFA is smaller than the PMI, $\alpha_i < \beta_{ji}$, $\alpha< \bar{\beta}_j$, which is typical in many applications, then asymptotics \eqref{FOAOmomentsgenCh8} and \eqref{FOAOmomentsgenCh8_2} are reduced to
\begin{equation}\label{FOAOmomentsgenCh8_3} 
 \inf_{\delta \in \class_\pi(\alphab,\betab)} \Rca^m_{i,\theta_i}(\delta)   \sim \brc{\frac{|\log \alpha_{i}|}{I_i(\theta_i)}}^m \sim \Rca^m_{i,\theta_i}(\delta_A) 
 \end{equation}
(as $\alpha_{\max}, \beta_{\max} \to 0$) and
 \begin{equation}\label{FOAOmomentsgenCh8_4} 
 \inf_{\delta \in \class_\pi^\star(\alpha,\betab)} \Rca^m_{i,\theta_i}(\delta)   \sim \brc{\frac{|\log \alpha|}{I_i(\theta_i)}}^m \sim \Rca^m_{i,\theta_i}(\delta_A) 
 \end{equation}
 (as $\alpha, \bar{\beta}_{\max} \to 0$).
 
Consider now the fully Bayesian setting where not only the prior distribution $\pi=\{\pi_k\}_{k\in\Zbb}$ of the changepoint $\nu$ is given, but also the prior distribution $p=\{p_i\}_{i\in\Nc}$ of hypotheses 
$\Pb(\Hyp_i)=p_i$, $i\in \Nc$ is specified. Then in place of the maximal probabilities of misidentification \eqref{PMIsupdef} one can consider the following average probabilities of misidentification
\begin{align*}
\PMI_i^{\pi,W}(\delta)& = \Pb^{\pi,W}_i(d \neq i, T<\infty| T>\nu)
\\
&= \int_{\Theta_i} \Pb_{i,\theta_i}^\pi(d \neq i, T<\infty| T>\nu) \mrm{d} W_i(\theta_i),
\\
\overline{\PMI}^{\pi,W,p}(\delta) & =\sum_{i=1}^N p_i \PMI_i^{\pi,W}(\delta),
\end{align*}
and the risk associated with the detection delay is measured by $\Rca_{\pi,W,p}^r(\delta) = \Eb^{\pi,W,p}[(T-\nu)^r | T>\nu]$ (in place of \eqref{Ch8Riskdefmultiple}). Here 
\begin{align*}
\Pb^{\pi,W}_i(\Ac\times\Kc) &= \sum_{k\in \Kc} \pi_k \int_{\Theta_i}  \Pb_{k, i,\theta_i}(\Ac) \mrm{d} W_i(\theta_i),
\\
\Pb^{\pi,W,p}(\Ac\times\Kc) &= \sum_{i=1}^N p_i \Pb^{\pi,W}_i(\Ac\times\Kc) ,
\end{align*}
and $\Eb^{\pi,W,p}$ is the expectation under the measure $\Pb^{\pi,W,p}$. It follows from Theorem~\ref{Th:UpperboundsPFAPMI} that for the rule $\delta_A$ with $A_{i0}=A_0$, $A_{ij}=A_j$, $i\in\Nc$, $j\in \Nc_0$
we have
\begin{align*}
\PMI_i^{\pi,W,p}(\delta)& \le \frac{1+A_{0}}{A_{0}} \frac{N-1}{A_{i}}, \quad i \in \Nc, \nonumber
\\ 
\overline{\PMI}^{\pi,W,p}(\delta)& \le \frac{(1+A_{0})(N-1)}{A_{0}} \sum_{i=1}^N \frac{p_i}{A_{i}}. 
\end{align*}
Introduce the class of detection--identification rules
\[
\bar{\class}_{\pi,W,p}(\alpha,\beta)=\set{\delta: \PFA^\pi(\delta)\le \alpha ~ \text{and}~\overline{\PMI}^{\pi,W,p}(\delta) \le \beta}
\]
for which the weighted probability of false alarm does not exceed $\alpha\in(0,1)$ and the average probability of misidentification does not exceed $\beta\in(0,1)$.  Note that 
$\delta_A\in \bar{\class}_{\pi,W,p}(\alpha,\beta)$ whenever 
\[
 A_{0} = \frac{N}{\alpha} (1-\alpha/N) ~~ \text{and}~~ A_{ij}=A_1 =\frac{N-1}{(1-\alpha/N)\beta}  .
\]
Using Theorem~\ref{Th:FOAOCh8} it is easy to prove that rule $\delta_A$ is first-order asymptotically optimal in the fully Bayesian setting in class $\bar{\class}_{\pi,W,p}(\alpha,\beta)$. 
Specifically, the following theorem holds.

%%Theorem
\begin{theorem}\label{Th:FOAOBayes} 
Let $r\ge 1$, let the prior distribution of the change point belong to class $\Cb(\mu)$, and let $p=\{p_i\}_{i\in\Nc}$ be the prior distribution of hypotheses that the change occurs in the $i$th data stream.  
Assume that for some  $0< I_i(\theta_i)<\infty$, $\theta_i\in\Theta_i$, $i\in \Nc$ and $0<I_{ij}(\theta_i,\theta_j)<\infty$, 
$\theta_i\in \Theta_i$, $\theta_j\in\Theta_j$, $i\in \Nc$, $j\in\Nc\setminus\{i\}$  right-tail and left-tail conditions
$\C_{1}$ and $\C_{2}$ are satisfied and that $\inf_{\theta_j \in \Theta_j} I_{ij}(\theta_i, \theta_j) > 0$ for all $j\in \Nc\setminus\{i\}$, $i \in \Nc$. 
If thresholds $A_{i0}=A_0$ and $A_{ij}=A_1$, $j \in \Nc\setminus\{i\}$, $i \in \Nc$ in rule $\delta_A$ are so selected that $\PFA^\pi(\delta_A) \le \alpha$, $\overline{\PMI}^{\pi,W,p}(\delta) \le \beta$ and 
$\log A_{0} \sim |\log\alpha|$,  $\log A_{1} \sim |\log \beta|$ as  $\alpha,  \beta \to 0$, in particular as $A_{0}=N(1-\alpha/N)/\alpha$ and $ A_{1} =(N-1)[(1-\alpha/N)\beta]^{-1}$, then $\delta_{A}$ 
is first-order asymptotically  optimal as $\alpha, \beta\to 0$ in class $\bar{\class}_{\pi,W,p}(\alpha,\beta)$,  minimizing moments of the detection delay up to order $r$: for all $0<m \le r$,
\begin{equation}\label{FOAOBayes} 
\begin{split}
& \inf_{\delta \in \bar{\class}_{\pi,W,p}(\alpha,\beta)} \Rca_{\pi,W,p}^m(\delta)  \sim\Rca_{\pi,W,p}^m(\delta_A)    
\\
& \sim  \max\set{ \gamma_0(p,W,\mu) \, |\log \alpha| , \gamma_1(p,W) \, |\log \beta|}^m
 \\
 & \qquad \text{as}~ \alpha, \beta \to 0 ,
 \end{split}
 \end{equation}
 where
\begin{align*}
\gamma_0(p,W,\mu) &= \sum_{i=1}^N p_i \int_{\Theta_i}  \frac{1}{I_i(\theta_i)+\mu} \mrm{d} W_i(\theta_i),
\\
\gamma_1(p,W) &= \sum_{i=1}^N p_i \int_{\Theta_i} \frac{1}{\displaystyle{\min_{{j\in\Nc\setminus\{i\}}}} \inf_{\theta_j\in\Theta_j} I_{ij}(\theta_i, \theta_j)} \mrm{d} W_i(\theta_i).
\end{align*}
\end{theorem}

%%%%%Remark
\begin{remark} \label{Rem:AccADD}
First-order approximations \eqref{FOAOmomentsgenkCh8}--\eqref{FOAOmomentsgenCh8_2} and \eqref{FOAOBayes} for the moments of the detection delay are usually not accurate. In the general non-i.i.d.\ case, 
it is difficult if at all possible, to obtain more accurate higher-order approximations. Higher-order approximations for the expected detection delays ($m=1$) can be obtained in the i.i.d.\ case using nonlinear 
renewal theory and techniques developed in  \cite[Th 3.3 ]{dtv2} and \cite[Th 4.3.4, Th 7.1.5]{TNB_book2014}.
\end{remark}

%%__________________________
\section{An Example: Detection of Signals with Unknown Amplitudes}\label{sec:Ex}

Suppose there is an $N$-channel sensor system and we are able to observe the output vector 
$X_n=(X_n(1),\dots, X_n(N))$, $n=1,2,\dots$ The observations $X_n(i)$ in the $i$th channel have the form
$$
X_{n}(i)=\theta_i S_{n}(i) \Ind{n > \nu}  +\xi_{n}(i),\quad n \ge 1, ~ i =1,\dots,N ,
$$
where $\theta_i$ is  an unknown intensity or amplitude ($\theta_i>0$)  of a deterministic signal $S_{n}(i)$ (e.g., the signal $S_{n}(i)= \cos (\omega_i n)$) and $\{\xi_{n}(i)\}_{n \in\Zbb_+}$, $i\in\Nc$ are 
mutually independent noises which are AR$(p)$ Gaussian stable processes that obey recursions
\begin{equation}\label{Ex1Ch8}
\xi_{n}(i) = \sum_{t=1}^p \varrho_{i,t} \xi_{n-t}(i) + w_{n}(i), \quad n \ge 1.
\end{equation}
Here $\{w_{n}(i)\}_{n\ge 1}$, $i\in \Nc$, are mutually independent i.i.d.\ Gaussian sequences with mean zero and standard deviation $\sigma>0$. 
The coefficients $\varrho_{i,1},\dots,\varrho_{i,p}$ and variance $\sigma^2$ are known.

A signal may appear only in one channel and should be detected and isolated quickly, i.e.,  the number of a channel where the signal 
appears should be identified along with detection.

Define $\widetilde{S}_{i,n} = S_{n}(i)- \sum_{t=1}^{p_n} \varrho_{i,t} S_{n-t}(i)$ and $\widetilde{X}_{i,n} = X_{n}(i) - \sum_{t=1}^{p_n} \varrho_{i,t} X_{n-t}(i)$ , where $p_n =p$ if $n > p$ and $p_n =n$ if $n \le p$.  
The LLRs have the form
\[
\begin{split}
\lambda_{i,\theta_i}(k, k+n) &= \frac{\theta_i}{\sigma^2}  \sum_{t=k+1}^{k+n} \wtS_{i,t} \wtX_{i,t} -\frac{\theta_i^2}{2 \sigma^2} \sum_{t=k+1}^{k+n} \wtS_{i,t}^2 ,
\\
\lambda_{i,\theta_i; j, \theta_j}(k, k+n) &= \lambda_{i,\theta_i}(k, k+n) - \lambda_{j,\theta_j}(k, k+n) .
\end{split}
\]
Under measure $\Pb_{k,i,\vartheta}$, $\vartheta \in \Theta_i$, the LLR $\lambda_{i,\theta_i; j, \theta_j}(k, k+n) $ is a Gaussian process (with independent non-identically distributed increments) 
with mean and variance
\begin{equation}\label{LLRARCh8}
\begin{split}
&\Eb_{k,i,\vartheta} [\lambda_{i, \theta_i; j, \theta_j}(k, k+n)]  = 
\\
& \frac{1}{2\sigma^2} \brcs{(2 \theta_i \vartheta -\theta_i^2)\sum_{t=k+1}^{k+n} \wtS_{i,t}^2 + \theta_j^2 \sum_{t=k+1}^{k+n} \wtS_{j,t}^2} ,
\\
& \Var_{k,i,\vartheta} [\lambda_{i, \theta_i; j, \theta_j}(k, k+n)]  =  
\\
& \frac{1}{\sigma^2} \brcs{\theta_i^2 \sum_{t=k+1}^{k+n} \wtS_{i,t}^2 +\theta_j^2 \sum_{t=k+1}^{k+n} \wtS_{j,t}^2} .
\end{split}
\end{equation}

Let $\Theta_i=(0,\infty)$, $i \in\Nc$ and assume that 
\begin{equation*}
\lim_{n\to \infty} \frac{1}{n}  \sup_{k \in \Zbb_+} \sum_{t=k+1}^{k+n} \wtS_{i,t}^2 = Q_i ,
\end{equation*}
where  $0<Q_i <\infty$. This is typically the case in most signal processing applications, e.g., for the sequence of sine pulses $S_{n}(i) = \sin(\omega_i n + \phi _i)$ with frequency $\omega_i$ and phase $\phi_i$. 
Then for all $k\in\Zbb_+$ and $\theta_i, \theta_j\in (0,\infty)$
\begin{align*}
&\frac{1}{n}\lambda_{i,\theta_i; j, \theta_j}(k, k+n)  \xra[n\to\infty]{ \Pb_{k,i,\theta_i} -\text{a.s.}} \frac{\theta_i^2 Q_i + \theta_j^2 Q_j}{2\sigma^2} 
\\
& =I_{ij}(\theta_i, \theta_j), \quad j \in \Nc\setminus\{i\}, ~ i \in \Nc,
\\
& \frac{1}{n}\lambda_{i,\theta_i}(k, k+n)  \xra[n\to\infty]{ \Pb_{k,i,\theta_i} -\text{a.s.}} \frac{\theta_i^2 Q_i}{2\sigma^2} 
\\
&=I_{i}(\theta_i), \quad  i \in \Nc,
\end{align*}
so that condition $\C_1$ holds. Furthermore, since all moments of the LLR are finite condition  $\C_2$ holds for all $r \ge 1$. Indeed, using \eqref{LLRARCh8}, we obtain that 
\[
I_i(\vartheta,\theta_i) =: \lim_{n\to\infty}  \frac{1}{n} \Eb_{k,i,\vartheta} [\lambda_{i,\theta_i}(k, k+n)]  = (\vartheta \theta_i - \vartheta^2/2) Q_i /\sigma^2
\]
and for any $\varkappa >0$
\begin{align*}
&\Pb_{k,i,\theta_i}\brc{\frac{1}{n} \inf_{|\vartheta-\theta_i| < \varkappa} \lambda_{i,\vartheta}(k, k+n)< I_i(\theta_i) - \varepsilon} \le
\\
&\Pb_{k,i,\theta_i}\brc{ \sup_{\vartheta \in [\theta_i -\varkappa, \theta_i+\varkappa]} \left\vert\frac{1}{n}\lambda_{i,\vartheta}(k, k+n)- I_i(\vartheta,\theta_i)\right\vert > \varepsilon} 
\\
&= \Pb_{k,i,\theta_i}\brc{|Y_{k, n}(\theta_i)| > \varepsilon \sqrt{n}},
\end{align*}
where
\[
Y_{k, n}(\theta_i) =  \frac{\theta_i}{\sqrt{n}} \sum_{t=k+1}^{k+n} \wtS_{i,t} \eta_t(i), \quad n \ge 1
\]
and $\{\eta_t(i)\}_{t \ge 1}$ is the sequence of standard zero-mean normal random variables. Hence $\{Y_{k, n}(\theta_i)\}_{n\ge 1}$ is the sequence of normal random variables with mean zero and 
variance $\sigma_{i,n}^2=  n^{-1}\sigma^{-2} \theta_i^2 \sum_{t=k+1}^{k+n} (\wtS_{i,t})^2$,  which is asymptotic to $\theta_i^2 Q_i/\sigma^2$. 
Thus, for a sufficiently large $n$ there exists $\delta_0>0$ such that  $\sigma^2_n \le \delta_0 + \theta_i^2 Q_i/\sigma^2$, and we obtain that for all large $n$
\begin{align*}
&\Pb_{k,i,\theta_i}\brc{\frac{1}{n} \inf_{|\vartheta-\theta_i| < \varkappa} \lambda_{i,\vartheta}(k, k+n)< I_i(\theta_i) - \varepsilon} \le
\\
&  \Pb\brc{|\hat\eta| > \frac{ \delta_0 + \theta_i^2 Q_i/\sigma^2}{\sigma_n^2} \frac{\varepsilon \sqrt{n}}{\delta_0 + \theta_i^2 Q_i/\sigma^2}} 
 \\
 & \le  \Pb\brc{|\hat\eta| > \frac{\varepsilon \sqrt{n}}{\delta_0 + \theta_i^2 Q_i/\sigma^2}},
\end{align*}
where $\hat\eta$ is a standard normal random variable. Therefore, 
\begin{align*}
& \Upsilon_r(\varkappa,\varepsilon;i,\theta_i) 
\\ 
& = \sum_{n=1}^\infty n^{r-1} \sup_{k\in\Zbb_+} \Pb_{k,i,\theta_i}\Bigg\{\frac{1}{n} \inf_{|\vartheta-\theta_i| < \varkappa} \lambda_{i,\vartheta}(k, k+n)
\\
& \quad \quad < I_i(\theta_i) - \varepsilon \Bigg\}
\\
& \le \sum_{n=1}^\infty n^{r-1}  \Pb\brc{|\hat\eta| > \frac{\varepsilon \sqrt{n}}{\delta_0 + \theta_i^2 Q_i/\sigma^2}},
\end{align*}
where the right-hand side term is finite for all $r\ge 1$ due to the finiteness of all moments of the normal distribution,
so that condition $\C_2$ holds for all $r \ge 1$.

Obviously, $\inf_{\theta_j\in (0,\infty)} I_{ij}(\theta_i,\theta_j) = \theta_i^2Q_i/(2\sigma^2)=I_i(\theta_i) > 0$.
Therefore, by Theorem~\ref{Th:FOAOCh8}, the detection--identification rule $\delta_A$ is asymptotically first-order optimal with respect to all positive moments of the detection delay and asymptotic formulas
\eqref{FOAOmomentsgenCh8} and \eqref{FOAOmomentsgenCh8_2} hold with 
\[
\inf_{\theta_j\in (0,\infty)} I_{ij}(\theta_i,\theta_j) = I_i(\theta_i) = \frac{\theta_i^2Q_i}{2\sigma^2}.
\] 
If $\max_{j \neq i} \beta_{ji} \ge \alpha_i$,  $\max_{j \neq i} \bar\beta_{j} \ge \alpha$, and $\mu=0$, then asymptotic formulas \eqref{FOAOmomentsgenCh8_3} and \eqref{FOAOmomentsgenCh8_4} hold.

Note that by condition $\C_2$ rule $\delta_A$ is asymptotically optimal for almost arbitrary mixing distributions $W_i(\theta_i)$. In this example, it is most convenient to 
select the conjugate prior, $W_i(\theta_i) =  F(\theta_i/v_i)$, where $F(y)$ is a standard normal distribution and $v_i>0$, in which case the decision statistics can be computed explicitly.

It is worth noting that this example arises in certain interesting practical applications, e.g., in multichannel/multisensor surveillance systems such as radars, sonars, and electro-optic/infrared sensor systems, which  deal 
with detecting moving and maneuvering targets that appear at unknown times, and  it is necessary to detect a signal from
a randomly appearing target in clutter and noise with the minimal average detection delay as well as to identify a channel where it appears. 
See \cite{Bakutetal-book63,Richards_radar2014,Marage_sonar2013,Tartakovsky&Brown-IEEEAES08}. 
Another challenging application area where the multichannel model is useful is cyber-security
\cite{Tartakovsky-Cybersecurity14,Tartakovskyetal-SM06,Tartakovskyetal_IEEESP2013}. Malicious intrusion attempts in computer networks (spam campaigns, personal data theft, 
worms, distributed denial-of-service (DDoS) attacks, etc.)  incur significant financial damage and are severe harm to the integrity of personal information.  
It is therefore essential to devise automated techniques to detect computer  network intrusions as quickly as possible so that an appropriate response can be provided and the 
negative consequences for the users are eliminated. In particular, DDoS attacks typically involve many traffic streams resulting in a large number of packets aimed at 
congesting the target's server or network.  

%%__________________________________________
\section{Concluding Remarks}\label{sec:Remarks}

%%1
1. Since we do not specify a class of models for the observations such as Gaussian, Markov, or HMM and build the decision statistics on the LLR processes, we restrict the behavior of LLRs which is expressed by conditions 
$\C_1$ and $\C_2$ related to the law of large numbers for the LLR and rates of convergence in the law of large numbers. As the example in Section~\ref{sec:Ex} shows, 
these conditions hold for the additive changes (in the mean) of the AR($p$) process governed by 
the Gaussian process. These conditions also hold in a variety of non-additive examples (detection of changes in spectrum of time series such as AR($p$) and ARCH($p$) processes) as well as for a large class of
homogeneous Markov processes~\cite{PergTarSISP2018,PergTarJMA2019},~\cite[Sec 3.1, Ch 4]{AGT_book2020} and for hidden Markov models with finite hidden state space~\cite{Fuh&TartakovskyIEEEIT2019}.

%%2
2. While we focused on the multistream detection--identification problem \eqref{modeldetisol}, it should be noted that similar results also hold in the ``scalar'' detection--isolation problem when the 
observations $\{X_n\}_{n \ge 1}$  represent either a scalar process or a vector process but 
all components of this process change at time $\nu$. Specifically, let $\{f_{\theta}(X_t|\Xb^{t-1}), \theta\in\Theta\}$ be a parametric family of densities and for  $i=1,\dots,N$ and $\Theta_i\subset\Theta$ consider the model 
\begin{equation*}
\begin{split}
& p(\Xb^n | \Hyp_{\nu,i}, \theta)  =  p(\Xb^n| \Hyp_\infty)= \prod_{t=1}^n f_{\theta_0}(X_t|\Xb^{t-1}) ~ \text{for}~ \nu \ge n ,
\\
& p(\Xb^n | \Hyp_{\nu,i},\theta)  =  \prod_{t=1}^{\nu} f_{\theta_0}(X_t|\Xb^{t-1})  \prod_{t=\nu+1}^{n}  f_{\theta}(X_t|\Xb^{t-1}) 
\\
 &\quad \text{for}~ \nu < n,~ \theta\in \Theta_i ,
\end{split}
\end{equation*}
where  $\theta_0$ is the known pre-change parameter and $\theta$ is the unknown post-change parameter. In other words,  there are $N$ types of change and for the $i$th type of change the 
value of the post-change parameter $\theta$ belongs to a subset $\Theta_i$ of the parameter space $\Theta$.  It is necessary to detect and isolate a change as rapidly as possible, i.e., 
to identify what type of change has occurred.  The change detection--identification rule $\delta_A=(d_A, T_A)$ is defined as in \eqref{CPDisolRule} where the statistics
$\bar{\Lambda}_{i j}^{\pi,W}(n)$ get modified as follows
\begin{equation*}
\begin{split}
\bar{\Lambda}_{i j}^{\pi,W}(n)& = \frac{\sum_{k=-1}^{n-1}\pi_k\int_{\Theta_i} LR_{\theta}(k,n) \, \D W_i(\theta)}{\sum_{k=-1}^{n-1}\pi_k \sup_{\theta\in\Theta_j} LR_{\theta}(k,n)}, ~~ i,j \in \Nc;
\\
\bar{\Lambda}_{i 0}^{\pi,W}(n)& = \frac{\sum_{k=-1}^{n-1}\pi_k\int_{\Theta_i} LR_{\theta}(k,n) \, \D W_i(\theta)}{\Pb(\nu \ge n)}, \quad i \in \Nc
\end{split}
\end{equation*}
with the likelihood ratio
\[
 LR_{\theta}(k,n) = \prod_{t=k+1}^{n} \frac{f_{\theta}(X_t|\Xb^{t-1})}{f_{\theta_0}(X_t|\Xb^{t-1})} .
\]

Write $\lambda_\theta(k,k+n) = \log  LR_{\theta}(k,n)$ and $\lambda_{\theta, \theta^*}(k,k+n)=  \lambda_{\theta}(k,k+n) - \lambda_{\theta^{*}}(k,k+n)$, where $ \lambda_{\theta^{*}}(k,k+n)= 0$ for $\theta^*=\theta_0$. Conditions $\Cb_1$ and $\Cb_2$ also get modified as follows

\vspace{2mm}

\noindent $\C_{1}$. {\em  There exist positive and finite numbers $I(\theta,\theta_0)=I(\theta)$, $\theta\in \Theta_i$, $i \in \Nc$ and $I(\theta,\theta^*)$, $\theta^*\in \Theta_j$, $j \in \Nc\setminus\{i\}$, $\theta\in \Theta_i$, 
$i \in \Nc$, such that for any  $\varepsilon >0$ and all  $k\in \Zbb_+$, $\theta\in\Theta_i$, $\theta^*\in\Theta_j$,  $j \in \Nc_0\setminus\{i\}$, $i \in \Nc$}
\begin{equation*}
\lim_{M\to\infty} p_{M,k}(\varepsilon; \theta, \theta^*) =0 .
\end{equation*}

\noindent $\C_{2}$. {\em  For any $\varepsilon>0$ and some $r\ge 1$}
\begin{equation*}
\Upsilon_r(\varepsilon; \theta)  < \infty \quad \text{for all}~  \theta\in\Theta_i, ~ i\in \Nc  ,
\end{equation*}
where
\begin{align*}
 & p_{M,k}(\varepsilon; \theta; \theta^*)   =
 \\
 &\Pb_{k,\theta}\set{\frac{1}{M}\max_{1 \le n \le M} \lambda_{\theta, \theta^*}(k, k+n) \ge (1+\varepsilon) I(\theta, \theta^*)},
\\
& \Upsilon_r(\varepsilon; \theta)  =   \lim_{\varkappa\to0}\sum_{n=1}^\infty n^{r-1} 
\\
& \sup_{k\in \Zbb_+} \Pb_{k,\theta}\Bigg\{\frac{1}{n} \inf_{\{\vartheta\in\Theta\,:\,\vert \vartheta-\theta\vert<\varkappa\}}\lambda_{\vartheta}(k, k+n) 
 < I(\theta)  - \varepsilon\Bigg\} .
\end{align*}

Essentially the same argument shows that all previous results hold in this case too. In particular, the assertions of Theorem~\ref{Th:FOAOCh8} are correct: as $\alpha_{\max}, \beta_{\max} \to 0$ 
for all $\theta\in\Theta_i$ and all $i\in\Nc$
\begin{align*}
&  \inf_{\delta \in \class_\pi(\alphab,\betab)} \Rca^m_{\theta}(\delta)   \sim   \Rca^m_{\theta}(\delta_A)  \sim
\\
&\max\set{\frac{|\log \alpha_{i}|}{I(\theta)+\mu}, \max_{j\in\Nc\setminus\{i\}} \frac{|\log \beta_{ji}|}{\inf_{\theta^*\in\Theta_j} I(\theta, \theta^*)}}^m,
 \end{align*}
i.e., the detection--identification rule $\delta_A$ is asymptotically optimal to first order. 

Note also that, in general, these asymptotics are not reduced to \eqref{FOAOmomentsgenCh8_3} even when $\alpha_i=\beta_{ji}$. Everything depends on the configuration of the hypotheses. 

%%3

3.  All previous results can be easily generalized for the case when the change points are different for different streams, i.e., when $\nu=\nu_i$ with prior distributions $\pi_k^{(i)}= \Pb(\nu_i=k)$, assuming that condition \eqref{Prior} 
for $\pi_k=\pi_k^{(i)}$ holds with $\mu=\mu_i$, $i \in \Nc$. Then in relations \eqref{Psii}, \eqref{Psii2}, \eqref{PsiACh8}, \eqref{FOAOmomentsgenkCh8}, \eqref{FOAOmomentsgenCh8}, 
\eqref{FOAOmomentsgenkCh8_2}, \eqref{FOAOmomentsgenCh8_2} and other relations where $\mu$ is present, the value of $\mu$ should be simply replaced with $\mu_i$.

%%4
4. For independent observations as well as for many Markov and certain hidden Markov models the decision statistics $\bar{\Lambda}_{i j}^{\pi,W}(n)$ can be computed effectively, 
so implementation of the proposed detection--identification rule is not
an issue. Still, in general, the computational complexity and memory requirements of rule $\delta_A$ are high. To avoid this complication,  
rule $\delta_A$ can be modified into a window-limited version where the summation in the statistics $\bar{\Lambda}_{i j}^{\pi,W}(n)$ over potential change points $k$ is restricted to the sliding window of size 
$\ell$.  Following guidelines of \cite[Ch 3, Sec 3.10]{AGT_book2020} (where asymptotic optimality of mixture window-limited rules was established in the single-stream case), 
it can be shown that the window-limited version
also has first-order asymptotic optimality properties as long as the size of the  window $\ell(A)$ goes to infinity as $A\to\infty$ at such a rate that $\ell(A)/\log A \to \infty$ but $\log \ell(A)/\log A \to 0$. The details are omitted.

%%4
5.  If $\pi\in \Cb(\mu=0)$ or $\pi^{\alpha,\beta}$ depends on $\alpha, \beta$ and $\mu_{\alpha,\beta}\to0$ as $\alpha_{\max},\beta_{\max}\to0$, then an alternative detection--identification rule 
$\delta^*_A=(d^*,T_A^*)$ defined as in \eqref{CPDisolRule}--\eqref{TAi} where in the  definition of $T_A^{(i)}$ the statistics $\bar{\Lambda}_{ij}^{\pi,W}(n)$ are replaced by the statistics
\begin{equation*}
\begin{split}
R_{ij}(n) & = \frac{\sum_{k=0}^{n-1} \int_{\Theta_i} LR_{i, \theta_i}(k,n) \, \D W_i(\theta_i)}{\sum_{k=0}^{n-1}\sup_{\theta_j\in\Theta_j} LR_{j,\theta_j}(k,n)}, ~ i,j \in\Nc;
\\
R_{i0}(n) & = \sum_{k=0}^{n-1}\int_{\Theta_i} LR_{i,\theta_i}(k,n) \, \D W_i(\theta_i), \quad i \in \Nc ,
\end{split}
\end{equation*}
is also asymptotically optimal to first order. Specifically, with a suitable selection of thresholds asymptotic approximations \eqref{FOAOmomentsgenCh8_3} and \eqref{FOAOmomentsgenCh8_4} hold for $\delta_A^*$.

%%5
6.  For practical purposes, it is more reasonable to consider a ``frequentist'' problem setup that does not use prior distributions of the changepoint $\pi$ and hypotheses $p$. We believe that the most reasonable
 performance metric for false alarms is the maximal conditional local probability of a false alarm in a prespecified time-window $\ell$, $\sup_{1 \le k <\infty}\Pb_\infty( k\le T< k+\ell| T >k) $ (see, e.g.,
 \cite{AGT_book2020,TNB_book2014} for a detailed discussion). The optimality results in the Bayesian problem obtained in this paper are of importance in the frequentist (minimax and pointwise) problem, which can 
 be embedded into the Bayesian criterion with asymptotically improper uniform distribution of the changepoint. See Pergamenchtchikov and Tartakovsky~\cite{PergTarSISP2018,PergTarJMA2019} 
 and Tartakovsky~\cite[Ch 4]{AGT_book2020} for the single population.

%%_________________________________________
\section*{Acknowledgement}

The author would like to thank referees whose comments improved the article.

%%_________________________________________________________
%\appendix

\renewcommand{\theequation}{A.\arabic{equation}}
\setcounter{equation}{0}

%%_________________________________________________________
\section*{Appendix: Proofs}

\begin{IEEEproof}[Proof of Theorem~\ref{Th:LBdetisol}]
The proof is split into two parts. 

%%%Part 1______________________________________________________________________________
{\em Part 1: Proof of asymptotic inequalities \eqref{LBinclassdetisolkPMI} and \eqref{LBinclassdetisolPMI}}.

To prove \eqref{LBinclassdetisolkPMI} and \eqref{LBinclassdetisolPMI} define 
\[
M_{\beta_{ji}}=M_{\beta_{ji}}(\varepsilon, \theta_i, \theta_j)= (1-\varepsilon) |\log\beta_{ji}|/I_{ij}(\theta_i,\theta_j)
\]
and note first that, by the  Chebyshev inequality, for every $\varepsilon \in(0,1)$ and $r >0$
\begin{align*}
& \Eb_{k, i,\theta_i}\left[(T-k)^r; d=i;  T> k \right] 
\\
 &\ge M_{\beta_{ji}}^r \Pb_{k,i,\theta_i}\set{T-k >  M_{\beta_{ji}}, d =i, T>k}
 \\
 & = M_{\beta_{ji}}^r \Pb_{k,i,\theta_i}\set{T-k >  M_{\beta_{ji}}, d =i} .
\end{align*}
Therefore,  for all $\theta_j\in\Theta_j$ and $j\in\Nc_0\setminus\{i\}$
\begin{align*}
& \Eb_{k,i,\theta_i}\left[(T-k)^r; d=i;  T> k \right]  
 \\
 &\ge \brcs{\frac{(1-\varepsilon)|\log \beta_{ji}|}{I_{ij}(\theta_i,\theta_j)}}^r (1+o(1))
\end{align*}
whenever for all $\varepsilon\in(0,1)$ and all fixed $k\in\Zbb$
\begin{equation}\label{Probkto1}
\lim_{\alpha_{\max}, \beta_{\max} \to 0} \inf_{\delta\in  \class_\pi(\alphab,\betab)} \Pb_{k, i,\theta_i}\set{T-k > M_{\beta_{ji}}, d =i} =1,
\end{equation}
and inequality \eqref{LBinclassdetisolkPMI} follows since $\varepsilon$ can be arbitrarily small and
\[
\Rc^r_{k,i,\theta_i}(\delta) \ge \Eb_{k,i,\theta_i}\left[(T-k)^r; d=i;  T> k \right].
\]
Recall that for $k=-1$ we set $T-k=T$ rather than $T+1$ everywhere. Note that $\Rc^r_{0,i,\theta_i}(\delta)\equiv \Rc^r_{-1,i,\theta_i}(\delta)$.

Analogously, 
\begin{align*}
& \Rca^r_{i,\theta_i}(\delta)  \ge \Eb_{i,\theta_i}^\pi \left[(T-\nu)^r; d=i;  T> \nu \right] 
\\
 &\ge  M_{\beta_{ji}}^r \Pb_{i,\theta_i}^\pi \set{T-\nu >  M_{\beta_{ji}}, d =i} ,
\end{align*}
so that inequality \eqref{LBinclassdetisolPMI} holds whenever 
\begin{equation}\label{Probpito1}
\lim_{\alpha_{\max}, \beta_{\max} \to 0} \inf_{\delta\in  \class_\pi(\alphab,\betab)} \Pb_{i,\theta_i}^\pi\set{T-\nu > M_{\beta_{ji}}, d =i} =1.
\end{equation}
Hence, we now focus on proving equalities \eqref{Probkto1} and  \eqref{Probpito1}.

Obviously, 
\begin{align*}
& \Pb_{k,i,\theta_i}(T-k>M_{\beta_{ji}}, d=i)  = \Pb_{k,i,\theta_i}(d =i)
\\
& \quad - \Pb_{k,i,\theta_i}(T-k\le M_{\beta_{ji}}, d=i)
\\
& =  1-\Pb_{k,i,\theta_i}(d \neq i) - \Pb_{i,k,\theta_i}(T \le k, d=i) 
\\
& \quad - \Pb_{k,i,\theta_i}(k < T\le M_{\beta_{ji}}+k, d=i),
\end{align*}
where $ \Pb_{i,k,\theta_i}(T \le k, d=i)= \Pb_{\infty}(T \le k, d=i)$.
Write $\Pi_k= \Pb(\nu >k)$. For any $\delta\in\class_\pi(\alphab,\betab)$ and $k\ge 0$, we have
\begin{align*}
\alpha_i & \ge \PFA_i(\delta) = \sum_{t=0}^\infty \pi_t \Pb_\infty(T \le t, d=i)  
\\
& \ge \sum_{t=k}^\infty \pi_t \Pb_\infty(T \le t, d=i)  \ge \Pb_\infty(T \le k, d=i) \Pi_{k-1},
\end{align*}
and
\[
\PMI_{ij}(\delta)= \sup_{\theta_i \in \Theta_i}\sum_{s=-1}^\infty \pi_s \Pb_{s,i,\theta_i}(d=j, T<\infty) \le \beta_{ij}
\]
so that, for any $\delta\in\class_\pi(\alphab,\betab)$, 
\begin{equation}\label{PTlessk}
 \Pb_\infty( T \le k, d=i) \le \alpha_i/\Pi_{k-1}, ~~ k \in \Zbb_+ 
\end{equation}
and
\begin{align}
\label{PMIkji}
\sup_{\theta_i \in \Theta_i} \Pb_{k,i,\theta_i}(d=j, T<\infty) \le  \beta_{ij}/\pi_k, ~~ k \in  \Zbb,
\\
\sup_{\theta_i \in \Theta_i} \Pb_{k,i,\theta_i}(d\neq i, T<\infty) \le  \pi_k^{-1} \sum_{j\in\Nc\setminus\{i\}} \beta_{ij} .
\label{PMIknoti}
\end{align}
Therefore,
\begin{align*}
& \Pb_{k,i,\theta_i}(T-k>M_{\beta_{ji}}, d=i)  
\\
& \ge  1- \alpha_i/\Pi_{k-1} -  \pi_k^{-1} \sum_{j\in\Nc\setminus\{i\}} \beta_{ij}
\\
& \quad - \Pb_{k,i,\theta_i}(k < T\le M_{\beta_{ji}}+k, d=i)  .
\end{align*}
This inequality implies that to prove  \eqref{Probkto1} we have to show that, as $\alpha_{\max}, \beta_{\max} \to 0$
\begin{equation}\label{Probkto0}
 \sup_{\delta \in \class_\pi(\alphab,\betab)} \Pb_{k,i,\theta_i}\set{0 <T-k \le  M_{\beta_{ji}}, d =i}  \to 0.
\end{equation}

 For the sake of brevity, we will write $\lambda_{i,j}(k,k+n)$ for the LLR $\lambda_{i,\theta_i; j, \theta_j}(k, k+n)$. Let $\Ac_{k,\beta} =\{k< T \le k + M_{\beta_{ji}}\}$ 
 and for $C>0$
 \[
 \Bc_{k,\beta,i,j} = \{d=i, \Ac_{k,\beta}\} \bigcap \{\max_{k < n \le k+M_{\beta_{ji}}}\lambda_{i,j}(k,n) <C\} .
 \]
 Changing the measure 
$\Pb_{k,j,\theta_j} \to  \Pb_{k,i,\theta_i}$, for any $C>0$ we obtain
\begin{align*}
& \Pb_{k,j,\theta_j}(d=i, T<\infty) = \Eb_{k,j,\theta_j}\brcs{\Ind{d=i,T<\infty}} 
\\
&= \Eb_{k,i,\theta_i}\brcs{\Ind{d=i,T<\infty} e^{-\lambda_{i,j}(k,T)}} 
\\
& \ge \Eb_{k,i,\theta_i}\brcs{\Ind{d=i, \Ac_{k,\beta}, \lambda_{i,j}(k,T) <C} e^{-\lambda_{i,j}(k,T)}} 
\\
&\ge  e^{-C} \, \Eb_{k,i,\theta_i}\brcs{\Ind{d=i, \Ac_{k,\beta}, \lambda_{i,j}(k,T) <C}} 
\\
&= e^{-C} \Pb_{k,i,\theta_i}\brc{\Bc_{k,\beta,i,j}}
\\
& \ge e^{-C} \, \Bigg[ \Pb_{k,i,\theta_i}(d=i, \Ac_{k,\beta}) 
\\
&\quad - \Pb_{k,i,\theta_i}\set{\max_{1\le n \le M_{\beta_{ji}}}\lambda_{i,j}(k,k+n) \ge C}\Bigg] ,
\end{align*}
where the last inequality follows from the trivial inequality $\Pb(A\cap B)\ge \Pb(A) - \Pb(B^c)$. It follows that
\begin{align*}
&\Pb_{k,i,\theta_i}(\Ac_{k,\beta}, d=i) \le \Pb_{k,j,\theta_j}(d=i, T<\infty)  e^{C}  
\\
& \quad + \Pb_{k,i,\theta_i}\set{\max_{1\le n \le M_{\beta_{ji}}}\lambda_{i,j}(k,k+n) \ge C}.
\end{align*}
Setting $C=(1+\varepsilon) I_{ij}(\theta_i,\theta_j)M_{\beta_{ji}}= (1-\varepsilon^2) |\log \beta_{ji}|$ yields
\begin{equation}\label{Probk}
\begin{split}
& \Pb_{k,i,\theta_i}\set{0 <T-k \le   M_{\beta_{ji}}, d =i}  
 \\
 &\le \Pb_{k,j,\theta_j}(d=i, T<\infty)  e^{(1-\varepsilon^2)  |\log \beta_{ji}|} 
 \\
 & \quad+ p_{M_{\beta_{ji}}, k} (\varepsilon; i, \theta_i; j,\theta_j) ,
 \end{split}
\end{equation}
where by \eqref{PMIkji}
\[
\sup_{\theta_j \in \Theta_j} \Pb_{k,j,\theta_j}(d=i, T<\infty) \le  \beta_{ji}/\pi_k,
\]
which along with \eqref{Probk} yields the inequality
\begin{equation*}%\label{Probk2}
\begin{split}
& \sup_{\delta \in \class_\pi(\alphab,\betab)} \Pb_{k,i,\theta_i}\set{0 <T-k \le M_{\beta_{ji}}, d =i}  \le 
\\
& \beta_{ji}^{\varepsilon^2}/\pi_k  +  p_{M_{\beta_{ji}}, k} (\varepsilon; i, \theta_i; j,\theta_j) .
\end{split}
\end{equation*}
The first term goes to zero for any fixed $k$ and the second term also goes to zero as $\beta_{\max}\to0$ by condition $\C_1$, which implies \eqref{Probkto0}
and \eqref{Probkto1}.

Next, multiplying both sides of inequality \eqref{Probk}  by $\pi_k$ and summing over $k\in \Zbb$, we obtain 
\begin{align*}
& \Pb_{i,\theta_i}^\pi\set{0 <T-\nu \le  M_{\beta_{ji}}, d =i}   \le \beta_{ji} e^{(1-\varepsilon^2)  |\log \beta_{ji}|} 
\\
& \quad + \sum_{k=-1}^{\infty} \pi_k p_{M_{\beta_{ji}}, k} (\varepsilon; i, \theta_i; j,\theta_j) 
\\
& \le \beta_{ji}^{\varepsilon^2} +  \Pb(\nu >K_{\beta}) 
+ \sum_{k=-1}^{K_{\beta}} \pi_k p_{M_{\beta_{ji}}, k} (\varepsilon; i, \theta_i; j,\theta_j) ,
\end{align*}
where $K_\beta$ is an arbitrary integer which goes to infinity as $\beta_{\max}\to 0$. 
Obviously, the first term goes to 0 as $\beta_{\max}\to 0$.  The second term $\Pb(\nu >K_{\beta})\to 0$ by conditions \eqref{Prior}  and \eqref{Prior1}.
The third term also goes to 0 due to condition $\Cb_1$ and Lebesgue's dominated convergence theorem. Hence, for any $\delta\in \class(\alphab,\betab)$,
\[
\Pb_{i,\theta_i}^\pi\set{0 <T-\nu \le  M_{\beta_{ji}}, d =i} \to 0 ~ \text{as} ~ \alpha_{\max}, \beta_{\max} \to 0 .
\]

Finally, we have 
\begin{align*}
& \Pb_{i,\theta_i}^\pi\set{T-\nu >  M_{\beta_{ji}}, d =i}  = \Pb_{i,\theta_i}^\pi(T >\nu, d=i) 
\\
& \quad - \Pb_{i,\theta_i}^\pi\set{0 <T-\nu \le  M_{\beta_{ji}}, d =i},
\end{align*}
where 
\begin{align}
& \Pb_{i,\theta_i}^\pi( T>\nu, d=i)  = \Pb_{i,\theta_i}^\pi( d=i |T>\nu) [1-\PFA^\pi(\delta)] \nonumber
\\
& \ge \brc{1-\sum_{j \in \Nc\setminus\{i\}} \beta_{ij}}\brc{1-\sum_{\ell=1}^N \alpha_\ell} \to 1
\label{LBTnudi}
\end{align}
as $\alpha_{\max},\beta_{\max}\to0$ for any $\delta\in\class_\pi(\alphab,\betab)$. This yields \eqref{Probpito1}, and therefore, inequalities  \eqref{LBinclassdetisolPMI}.

%%---------------------------------------------------------------------------------------------------------------------------------------------
{\em Part 2: Proof of asymptotic inequalities  \eqref{LBinclassdetisolkPFA} and \eqref{LBinclassdetisolPFA}}.

Changing the measure $\Pb_\infty\to \Pb_{k,i,\theta_i}$ and using an argument similar to that used in Part 1 to obtain \eqref{Probk} with $M_{\beta_{ij}}$ replaced by
\[
N_{\alpha_i} = \frac{(1-\varepsilon) |\log \alpha_i|}{I_i(\theta_i) + \mu + \varepsilon_1}
\]
 we obtain
 \begin{align} \label{Pkiineq}
& \Pb_{k, i,\theta_i}\set{0 <T-k \le N_{\alpha_i}, d =i}   \nonumber
\\
& \le  e^{(1+\varepsilon)  I_i(\theta_i)  N_{\alpha_i} } \Pb_{\infty}\set{0 <T-k \le N_{\alpha_i}, d =i} \nonumber
\\
& +  \Pb_{k,i,\theta_i}\set{\frac{1}{N_{\alpha_i}}\max_{1\le n \le N_{\alpha_i}}\lambda_{i,\theta_i}(k,k+n) \ge (1+\varepsilon) I_i(\theta_i)} ,
\end{align}
where for all $\varepsilon \in (0,1)$
\begin{align} \label{IneqU1Ch8}
& e^{(1+\varepsilon)  I_i(\theta_i)  N_{\alpha_i} } \Pb_{\infty}\set{0 <T-k \le N_{\alpha_i}, d =i}  \nonumber
\\
& \le \exp\set{- \varepsilon^2 |\log \alpha_i| + (\mu+\varepsilon_1) (k-1)}   \nonumber
\\
&:= \overline{U}_{\alpha_i,k}(\varepsilon,\varepsilon_1) .
\end{align}

Using \eqref{Pkiineq} and \eqref{IneqU1Ch8}, we obtain
\begin{align*}
& \sup_{\delta\in\class_\pi(\alphab,\betab)} \Pb_{k, i,\theta_i}\set{0 <T-k \le N_{\alpha_i}, d =i} \le \overline{U}_{\alpha_i,k}(\varepsilon,\varepsilon_1) 
\\
&\qquad + p_{N_{\alpha_i},k}(\varepsilon;i,\theta_i),
\end{align*}
where for every fixed $k\in\Zbb_+$ the value of $ \overline{U}_{\alpha_i,k}(\varepsilon,\varepsilon_1)$ tends to zero and also $ p_{N_{\alpha_i},k}(\varepsilon;i,\theta_i)\to 0$ as $\alpha_{\max}\to0$ by
condition $\C_1$. Hence, it follows that for every fixed $k\in\Zbb$
\begin{equation}\label{supProbkto0}
\lim_{\alpha_{\max}\to0} \sup_{\delta\in\class_\pi(\alphab,\betab)} \Pb_{k, i,\theta_i}\set{0 <T-k \le N_{\alpha_i}, d =i}  =0.
\end{equation}

Next, we have
\begin{align*}
\Rc_{k, i, \theta_i}^r(\delta) &\ge \Eb_{k, i, \theta_i}[( T-k)^r, d=i, T> k] 
\\
&\ge N_{\alpha_i}^r \Pb_{k, i, \theta_i}(T  -k > N_{\alpha_i}, d=i, T>k)
\\
&= N_{\alpha_i}^r \Pb_{k, i, \theta_i}(T  -k > N_{\alpha_i}, d=i)
\\
& \ge  N_{\alpha_i}^r[\Pb_{k, i,\theta_i}( T>k, d=i)
\\
&\quad -\Pb_{k,i,\theta_i}(0 <  T -k \le  N_{\alpha_i}, d=i)] ,
\end{align*} 
where the second inequality follows from the Chebyshev inequality and 
\begin{align*}
& \Pb_{k, i,\theta_i}( T>k, d=i)=1-\Pb_{k,i,\theta_i}(d \neq i) 
\\
& \quad - \Pb_{i,k,\theta_i}(T \le k, d=i)  
\\
&\ge 1- \alpha_i/\Pi_{k-1} -  \pi_k^{-1} \sum_{j\in\Nc\setminus\{i\}} \beta_{ij}
\end{align*}
(see \eqref{PMIkji}--\eqref{PMIknoti}). Therefore,
\begin{align*}
& \inf_{\delta\in\class_\pi(\alphab,\betab)}\Rc_{k, i, \theta_i}^r(\delta) \ge 
 N_{\alpha_i}^r\big[ \inf_{\delta\in\class_\pi(\alphab,\betab)} \Pb_{\infty}( T>k, d=i)
\\
&\quad -\sup_{\delta\in\class_\pi(\alphab,\betab)}\Pb_{k,i,\theta_i}(0 <  T -k \le  N_{\alpha_i}, d=i)\big] ,
\end{align*}
where 
\begin{equation}\label{Pinfto0}
\lim_{\alpha_{\max}\to0} \inf_{\delta\in\class_\pi(\alphab,\betab)} \Pb_{\infty}( T>k, d=i) =1,
\end{equation}
and by \eqref{supProbkto0} the second term on the right hand-side goes to $0$ for any fixed $k \in\Zbb$.

It follows that for all fixed $k\in\Zbb$
\[
 \inf_{\delta\in\class_\pi(\alphab,\betab)} \Rc_{k, i,\theta_i}^r(\delta)  \ge \brcs{\frac{(1-\varepsilon)|\log\alpha_i|}{I_i(\theta_i)+\mu+\varepsilon_1}}^r (1+o(1)),
\]
where $\varepsilon$ and $\varepsilon_1$ can be arbitrarily small, which implies the inequality \eqref{LBinclassdetisolkPFA}.

Next, define
\[
K_{\alpha_i}=K_{\alpha_i}(\varepsilon,\mu,\varepsilon_1)=  \left\lfloor  \frac{\varepsilon^3 |\log\alpha_i|}{\mu+\varepsilon_1}\right\rfloor .
\]
Using inequalities \eqref{Pkiineq} and \eqref{IneqU1Ch8}, we obtain 
\begin{align*}
& \Pb^\pi_{i,\theta_i}(0 <  T -\nu \le  N_{\alpha_i}, d=i)  
\\
&= \sum_{k=-1}^\infty  \pi_k   \Pb_{k,i,\theta_i}\brc{0<  T -k \le   N_{\alpha_i}, d=i} \nonumber
\\
&= \sum_{k=-1}^{K_{\alpha_i}}  \pi_k   \Pb_{k,i,\theta_i}\brc{0<  T -k \le  N_{\alpha_i}, d=i} 
\\& \quad + 
\sum_{k=K_{\alpha_i}+1}^\infty  \pi_k   \Pb_{k,i, \theta_i}\brc{0<  T -k\le  N_{\alpha_i}, d=i}  \nonumber
 \\
 & \le   \sum_{k=-1}^{K_{\alpha_i}} \pi_k \overline{U}_{\alpha_i,k}(\varepsilon,\varepsilon_1)+  \sum_{k=-1}^{K_{\alpha_i}}  \pi_k  p_{N_{\alpha_i}, k} (\varepsilon; i, \theta_i)
 \\
 & \qquad + \sum_{k=K_{\alpha_i}+1}^\infty  \pi_k 
 \\
 &\le \Pi_{K_{\alpha_i}} + \max_{-1 \le k \le K_{\alpha_i}}  \overline{U}_{\alpha_i,k}(\varepsilon,\varepsilon_1) + \sum_{k=-1}^{K_{\alpha_i}}  \pi_k  p_{N_{\alpha_i}, k} (\varepsilon; i, \theta_i) \nonumber
 \\
 & = \Pi_{K_{\alpha_i}} +  \overline{U}_{\alpha_i,K_{\alpha_i}}(\varepsilon,\varepsilon_1) +  \sum_{k=-1}^{K_{\alpha_i}}  \pi_k  p_{N_{\alpha_i}, k} (\varepsilon; i, \theta_i),
 \end{align*}
 where $T-k=T$ for $k=-1$.
If $\mu>0$, by condition \eqref{Prior},  $\log \Pi_{K_{\alpha_i}} \sim - \mu \, K_{\alpha_i}$ as $\alpha_{\max} \to0$, so $\Pi_{K_{\alpha_i}}\to 0$. If $\mu=0$,  
this probability goes to $0$ as $\alpha_{\max} \to 0$ as well since, by condition \eqref{Prior1}, 
\[
\Pi_{K_{\alpha_i}} < \sum_{k= K_{\alpha_i}}^\infty  \pi_k |\log \pi_k| \xra[\alpha_{\max}\to0]{} 0. 
\]
Obviously, the second term $\overline{U}_{\alpha_i,K_{\alpha_i}}(\varepsilon,\varepsilon_1)\to0$ as $\alpha_{\max}\to0$. By condition $\C_1$ and Lebesgue's dominated convergence theorem, the third term goes 
to 0, and therefore,  all three terms go to zero as $\alpha_{\max}, \beta_{\max}\to0$ for all $\varepsilon,\varepsilon_1>0$, so that 
 \[
 \Pb^\pi_{i,\theta_i}(0 <  T -\nu \le  N_{\alpha_i}, d=i)  \to 0 \quad \text{as}~ \alpha_{\max}, \beta_{\max} \to 0.
 \]  
Since
\begin{align*}
&\Pb_{i, \theta_i}^\pi(T  -\nu > N_{\alpha_i}, d=i) 
\\
&= \Pb_{i,\theta_i}^\pi( T>\nu, d=i)-\Pb_{i,\theta_i}^\pi(0 <  T-\nu \le N_{\alpha_i}, d=i)
\end{align*}
and by \eqref{LBTnudi} $\Pb_{i,\theta_i}^\pi( T>\nu, d=i)   \to 1$ as $\alpha_{\max},\beta_{\max}\to0$ for any $\delta\in\class_\pi(\alphab,\betab)$, it follows that
\[
\Pb_{i, \theta_i}^\pi(T  -\nu > N_{\alpha_i}, d=i) \to 1 \quad \text{as}~ \alpha_{\max}, \beta_{\max} \to 0.
\]
Finally, by the Chebyshev inequality,
\begin{align*}
\Rca_{i, \theta_i}^r(\delta) &\ge \Eb_{i, \theta_i}^\pi[( T-\nu)^r, d=i, T> \nu] 
\\
&\ge N_{\alpha_i}^r \Pb_{i, \theta_i}^\pi(T  -\nu > N_{\alpha_i}, d=i),
\end{align*}
which implies that for any $\delta\in\class_\pi(\alphab,\betab)$ as $\alpha_{\max}, \beta_{\max} \to 0$
\[
\Rca_{i,\theta_i}^r(\delta)  \ge \brcs{\frac{(1-\varepsilon)|\log\alpha_i|}{I_i(\theta_i)+\mu+\varepsilon_1}}^r (1+o(1)) .
\]
Owing to the fact that $\varepsilon$ and $\varepsilon_1$ can be arbitrarily small the inequality \eqref{LBinclassdetisolPFA} follows.
\end{IEEEproof}

%%Proof of Lemma 1
\begin{IEEEproof}[Proof of Lemma~\ref{Lem:UpperEki}]
For $k \in \Zbb_+$, define the exit times
\begin{align*}
\tau_{i}^{(k)}(A)  & =  \inf \{n\ge 1: \lambda_{i,W}(k,k+n) - \lambda_{j}^\pi(k+n) \ge 
\\
&\log (A_{ij}/\pi_k) ~ \forall ~ j \in \Nc_0\setminus\{i\}\}, ~ i \in \Nc,
\end{align*}
where $\lambda_{i,W}(k,k+n)=\log \Lambda_{i,W}(k,k+n)$ and $\lambda_{0}^\pi(k+n) =\log \Pb(\nu \ge k+n)= \log \Pi_{k+n-1}$.

Obviously, for any $n >k$ and $k\in\Zbb_+$,
\begin{align*}
\log \bar{\Lambda}_{ij}^{\pi,W}(n)   & \ge \log\brc{\frac{\pi_k  LR_{i,W}(k,n)}{\sum_{\ell=-1}^{n-1}\pi_\ell \sup_{\theta_j\in\Theta_j} LR_{j,\theta_j}(\ell,n)}} 
\\
& = \lambda_{i,W}(k,n)  -\lambda_{j}^\pi(n) + \log \pi_k,
\end{align*}
so  for every set $A=(A_{ij})$ of positive thresholds $A_{ij}$, we have $(T_A-k)^+ \le (T_A^{(i)}-k)^+  \le  \tau_i^{(k)}(A)$ and, hence, $\Eb_{k,i,\theta_i}[(T_A-k)^+]^r\le \Eb_{k,i,\theta_i} [(\tau_i^{(k)}(A))^r]$.
Note that since we set $T_A-k=T_A$ for $k=-1$, it follows that $\Eb_{-1,i,\theta_i}[(T_A-k)^+]^r=\Eb_{0,i,\theta_i}[T_A]^r  \le  \Eb_{0,i,\theta_i} [(\tau_i^{(0)}(A))^r]$.

Setting $\tau=\tau_i^{(k)}(A)$ and $N=M_i(A)$ in inequality~(A.1) in Lemma~A1 in \cite[p. 239]{AGT_book2020}
we obtain that the following inequality holds: 
\begin{equation}\label{EktauineqCh8}
\begin{split}
& \Eb_{k, i, \theta_i} \brcs{\brc{\tau_i^{(k)}(A)}^r}   \le [M_i(A)]^{r} 
\\
&\quad + r 2^{r-1} \sum_{n=M_i(A)}^{\infty}  n^{r-1}   \Pb_{k, i, \theta_i}\brc{\tau_i^{(k)}(A) >  n}.
 \end{split}
 \end{equation}
 Next, we have
\begin{align*}
&\Pb_{k, i,\theta_i}\brc{ \tau_i^{(k)}(A) >n}  \le
\\
&\Pb_{k,i, \theta_i}\Bigg\{\frac{\lambda_{i,W}(k,k+n)-\lambda_{j}^\pi(k+n)}{n}  
\\
&\quad < \frac{1}{n} \log \brc{\frac{A_{ij}}{\pi_k}}, ~\forall ~j \in \Nc_0\setminus\{i\}\Bigg\} 
\\
& \le  \Pb_{k,i, \theta_i}\Bigg\{\frac{\lambda_{i,W}(k,k+n)- \log \Pi_{k+n-1}}{n} 
\\
& \quad < \frac{1}{n} \log \brc{\frac{A_{i0}}{\pi_k}}\Bigg\} .
\end{align*}
Let
\[
\widetilde{M}_i(A_{i0})=1+\left\lfloor \frac{\log (A_{i0}/\pi_k)}{I_i(\theta_i) + \mu - \varepsilon} \right \rfloor .
\]
Clearly, for all $n\ge \widetilde{M}_i(A_{i0})$ the last probability does not exceed the probability
\[
\Pb_{k,i, \theta_i}\set{\frac{\lambda_{i,W}(k,k+n)}{n} <I_i(\theta_i) + \mu -\varepsilon  -\frac{|\log\Pi_{k+n-1}|}{n}}
\]
and, by condition $\CP_1$, for a sufficiently large value of $A_{i0}$ there exists a small $\kappa$ such that 
\[
\left |\mu - \frac{|\log \Pi_{k+ \widetilde{M}_i(A_{i0})-1}|}{\widetilde{M}_i(A_{i0})} \right | < \kappa.
\]  
Therefore, for all sufficiently large $n$, 
\begin{align*}
&\Pb_{k,i, \theta_i}\brc{\tau_i^{(k)}(A) >n} 
\\
&\le \Pb_{k,i, \theta_i}\brc{\frac{1}{n}\lambda_{i,W}(k,k+n) <I_i(\theta_i) -\varepsilon + \kappa} .
\end{align*}
Also, 
\[
\lambda_{i,W}(k,k+n) \ge \inf_{\vartheta \in \Gamma_{\varkappa,\theta_i}} \lambda_{i,\vartheta}(k,k+n) +\log W_i(\Gamma_{\varkappa,\theta_i}),
\]
where $\Gamma_{\varkappa,\theta_i}=\{\vartheta\in\Theta_i\,:\,\vert\vartheta-\theta_i\vert<\varkappa\}$. Thus, for all sufficiently large  $n$ and $A_{\min}$, for which
$\kappa +|\log W(\Gamma_{\varkappa,\theta_i})|/n < \varepsilon/2$, we have
\begin{align}
&\Pb_{k,i, \theta_i}\brc{ \tau_i^{(k)}(A) >n}   \le \Pb_{k,i, \theta_i}\Bigg\{\frac{1}{n}\inf_{\vartheta\in\Gamma_{\varkappa,\theta_i}} \lambda_{i,\vartheta}(k,k+n) \nonumber
\\
&< I_i(\theta_i)  - \varepsilon+ \kappa + \frac{1}{n}|\log W(\Gamma_{\varkappa,\theta_i})|\Bigg\} \nonumber
\\
&  \le \Pb_{k,i, \theta_i}\brc{\frac{1}{n}\inf_{\vartheta\in\Gamma_{\varkappa,\theta_i}} \lambda_{i,\vartheta}(k,k+n) < I_i(\theta_i)  - \varepsilon/2}. \label{Probktau1Ch8}
\end{align}

Using \eqref{EktauineqCh8} and \eqref{Probktau1Ch8} yields inequality \eqref{EkineqCh8} and the proof is complete.
\end{IEEEproof}

%%Proof of Proposition
\begin{IEEEproof}[Proof of Proposition~\ref{Pr:FOapproxAdetisol}]
By Theorem~\ref{Th:UpperboundsPFAPMI}, the rule $\delta_A$ belongs to class $\class_\pi(\alphab,\betab)$ when 
\[
\alpha_i= \frac{1}{1+A_{i0}} ; \quad \beta_{ij} = \frac{1+A_{i0}}{A_{i0} \, A_{ji}}, ~~ j \in \Nc\setminus\{i\}, ~ i \in \Nc,
\]
and hence,  Theorem~\ref{Th:LBdetisol} implies (under condition $\Cb_1$) the asymptotic (as $A_{\min}\to\infty$) lower bounds
\begin{equation}\label{LBAdetisolk}
 \Rc^r_{k, i,\theta_i}(\delta_A) \ge [\Psi_i(A,\theta_i,\mu)]^r (1+o(1)) ~~ \forall~ k\in\Zbb
\end{equation}
and
\begin{equation}\label{LBAdetisol}
 \Rca^r_{i,\theta_i}(\delta_A) \ge [\Psi_i(A,\theta_i,\mu)]^r (1+o(1)) ,
\end{equation}
which hold for all  $ r>0$, $\theta_i\in \Theta_i$, and $i\in \Nc$. 
Thus, to prove the validity of the asymptotic approximations  \eqref{MADDAsaprk} and \eqref{MADDAsapr} it suffices to show that, under the left-tail condition $\C_2$, for $0 <m \le r$ and all 
$\theta_i\in\Theta_i$ and $i\in\Nc$ the following asymptotic upper bounds hold as $A_{\min}\to\infty$:
\begin{equation}\label{UppergenACh8k}
\Rc_{k, i, \theta_i}^m(\delta_A) \le   [\Psi_i(A,\theta_i,\mu)]^m (1+o(1)) ~~ \forall ~ k\in\Zbb
\end{equation}
and
\begin{equation}\label{UppergenACh8}
\Rca_{i, \theta_i}^m(\delta_A) \le   [\Psi_i(A,\theta_i,\mu)]^m (1+o(1)) . 
\end{equation}

It follows from inequality \eqref{EkineqCh8} in Lemma~\ref{Lem:UpperEki} that for any $0<\varepsilon < J_{ij}(\theta_i,\mu)$
\begin{align}\label{Ekiineq}
&  \Eb_{k,i,\theta_i}\brcs{\brc{T_A-k}^r; d_A=i; T_A> k}  \le  \nonumber
\\
& \brcs{1+\widetilde{\Psi}_i(A, \pi_k, \theta_i,\mu, \varepsilon)}^r  
 + r 2^{r-1}  \Upsilon_{r}(\varkappa, \varepsilon; i, \theta_i),
\end{align}
where $\Upsilon_{r}(\varkappa, \varepsilon; i, \theta_i)$ is defined in \eqref{UpsilonCh8}. Similarly to \eqref{PTlessk} we have $\Pb_\infty(T_A\le k, d_A=i) \le [(1+A_{i0})\Pi_{k-1}]^{-1}$, so that
\begin{align*}
\Pb_\infty(T_A\le k) \le \frac{1}{\Pi_{k-1}}\sum_{i=1}^N \frac{1}{1+A_{i0}},
\end{align*}
and hence,
\begin{align*}%\label{PinftyTAk}
\Pb_\infty(T_A > k) \ge 1- \frac{1}{\Pi_{k-1}}\sum_{i=1}^N \frac{1}{1+A_{i0}}.
\end{align*}
Using this inequality and inequality \eqref{Ekiineq}, we obtain
\begin{align}
 & \Rc_{k, i,\theta_i}^r(\delta_A) =\frac{\Eb_{k,i,\theta_i}\brcs{\brc{T_A-k}^r; d_A=i; T_A>k}}{\Pb_\infty(T_A > k)}  \nonumber
  \\
  &\le   \frac{\brc{1+\left\lfloor \frac{\log (A/\pi_k)}{I_i(\theta_i)+\mu-\varepsilon} \right\rfloor}^r + r 2^{r-1} \,\Upsilon_{r}(\varkappa, \varepsilon; i, \theta_i)}{1- \sum_{i=1}^N 1/[(1+A_{i0}) \Pi_{k-1}]}. \label{Rckupper}
 \end{align}
Since, by condition $\C_2$, $\Upsilon_{r}(\varkappa,\varepsilon; i, \theta_i) < \infty$ for all  $\theta_i\in\Theta_i$ and $i\in\Nc$, this implies the asymptotic upper bound \eqref{UppergenACh8k}. This completes the proof of the 
asymptotic approximation \eqref{MADDAsaprk}.

Next, using inequality \eqref{Ekiineq} we obtain
\begin{equation*}%\label{UpperExprCh8}
\begin{split}
&\Eb^\pi_{i,\theta_i} [(T_{A}-\nu)^r; d_A=i; T_A>\nu]  
\\
&= \sum_{k=-1}^\infty \pi_k  \Eb_{k, i, \theta_i}\brcs{(T_A-k)^r; d_A=i; T_A>k}
\\
&\le \sum_{k=-1}^\infty \pi_k \brcs{1+\widetilde{\Psi}_i(A, \pi_k, \theta_i,\mu, \varepsilon)}^r  
 + r 2^{r-1}  \Upsilon_{r}(\varkappa,\varepsilon; i, \theta_i).
\end{split}
\end{equation*}
 Recall that we set $T_A-k=T_A$ for $k=-1$. Applying this inequality together with inequality
 \[
1- \PFA^\pi(\delta_A) \ge 1- \sum_{i=1}^N \frac{1}{1+A_{i0}}
 \] 
(see \eqref{UpperPFAi})  yields
\begin{equation}\label{UpperRcaCh8}
\begin{split}
&\Rca_{i, \theta_i}^r(\delta_A) = \frac{\displaystyle{\sum_{k=-1}^\infty \pi_k \Eb_{k,i,\theta_i}[(T_A-k)^r; d_A=i;T_A>k]}}{1-\PFA^\pi(T_A)} 
\\
&\le\frac{\displaystyle{\sum_{k=-1}^\infty} \pi_k \brcs{1+\widetilde{\Psi}_i(A, \pi_k, \theta_i,\mu, \varepsilon)}^r  + r 2^{r-1}  \Upsilon_{r}(\varkappa,\varepsilon; i, \theta_i)}{1-\sum_{i=1}^N(1/(1+A_{i0})} .
\end{split}
\end{equation}
By condition $\C_2$, $\Upsilon_{r}(\varkappa,\varepsilon; i, \theta_i)< \infty$ for any $\varepsilon >0$ and any $\theta_i\in\Theta_i$ 
and, by condition \eqref{Prior1}, $\sum_{k=0}^\infty \pi_k |\log\pi_k|^r < \infty$. This implies that, as $A_{\min} \to \infty$, for all $0< m \le r$, all $\theta_i\in\Theta_i$, and all $i\in\Nc$ the following upper bound holds
\[
\Rca_{i, \theta_i}^r(\delta_A) \le \brcs{\widetilde{\Psi}_i(A, \pi_k=1, \theta_i,\mu, \varepsilon)}^r (1+o(1)) .
\]
Since $\varepsilon$ can be arbitrarily small and $\lim_{\varepsilon\to0} \widetilde{\Psi}_i(A, \pi_k=1, \theta_i,\mu, \varepsilon)=\Psi_i(A,\theta_i,\mu)$, the upper bound \eqref{UppergenACh8} follows and the proof of 
the asymptotic approximation  \eqref{MADDAsapr} is complete.
\end{IEEEproof}

%%%%%%%%%%%%%%%%%%%%%%%%%%%%%%%%%%%%%%%%%%%%%%
%\bibliographystyle{IEEEtranS}
%\bibliography{Bib_main}

% Generated by IEEEtranS.bst, version: 1.14 (2015/08/26)

%%%%% Biography%%%%%%%%
%\begin{IEEEbiography}%[{\includegraphics[width=1.2in,height=1.45in,keepaspectratio]{Tartakovsky-photo.pdf}}]
%{Alexander G.\ Tartakovsky} (M'01-SM'02)  

\begin{IEEEbiographynophoto}{Alexander G.\ Tartakovsky} (M'01-SM'02),  M.S., Ph.D., D.Sc., is an award-winning statistician, Head of the Space Informatics Laboratory at the Moscow Institute of Physics and Technology 
(``PhysTech''), and President of AGT StatConsult, Los Angeles, CA. From 2013 to 2015, he was a Professor of Statistics at the University of Connecticut, Storrs. Previously, for almost two decades, 
he was a Professor in the Department of Mathematics and the Associate Director of the Center for Applied Mathematical Sciences at the University of Southern California (USC).

Dr. Tartakovsky is the author of three books, several book chapters, and over 100 papers across a range of subjects, including theoretical and applied statistics; applied probability; sequential analysis; and changepoint detection. 
His research has many applications, including in statistical image and signal processing, video tracking, detection and tracking of targets in radar and infrared search and track systems, near-Earth space informatics, 
information integration/fusion, intrusion detection and network security, rapid detection of epidemics, and detection and tracking of malicious activity.

Dr. Tartakovsky earned an M.S. in Electrical Engineering from the Moscow Aviation Institute in 1978 and a Ph.D. in Statistics and Information Theory from PhysTech in 1981. 
He also earned an advanced Doctor of Science (D.Sc.) degree from PhysTech in 1990.

From 1981 to 1992, he was first a Senior Research Scientist and then Department Head at the Moscow Institute of Radio Technology and a Professor at PhysTech, 
where he worked on the application of statistical methods to optimization and modeling of information systems.

From 1993 to 1996, Dr. Tartakovsky was a professor at the University of California, Los Angeles (UCLA), first in the Department of Electrical Engineering and then in the Department of Mathematics.

Dr. Tartakovsky has received numerous awards for his work, including the Abraham Wald Prize in Sequential Analysis and several Best Young Scientist awards from the Russian Academy of Sciences. 
He is also a Fellow of the Institute of Mathematical Statistics (IMS) and a senior member of the Institute of Electrical and Electronics Engineers (IEEE).
\end{IEEEbiographynophoto}

\vfill

\end{document}